\documentclass{amsart}

\usepackage[T1]{fontenc}
\usepackage{graphicx}
\usepackage{mathptmx}  
\usepackage{amssymb}
\usepackage{amsmath}
\usepackage{cancel}
\usepackage{ulem}
\usepackage{color}
\usepackage{xcolor}
\usepackage{setspace}
\usepackage{nomencl}
\makenomenclature

\usepackage{kantlipsum}
\allowdisplaybreaks

\makeatletter
\def\subsection{\@startsection{subsection}{2}%
	\z@{.5\linespacing\@plus.7\linespacing}{.3\linespacing}%
	{\normalfont\bfseries}}
\makeatother

\newtheorem{theorem}[subsubsection]{Theorem}

\newtheorem{proposition}[subsubsection]{Proposition}
\newtheorem{definition}[subsubsection]{Definition}
\newtheorem{remark}[subsubsection]{Remark}
\newtheorem{example}[subsubsection]{Example}
\newtheorem{corollary}[subsubsection]{Corollary}

\title[Fractional finite sums theory: a motivation and some new results]{Fractional finite sums theory: \\ a motivation and some new results}

\author[L. F. Bielinski]{Leonardo F. Bielinski \\ {\tiny Universidade Estadual de Ponta Grossa \\ Ponta Grossa, PR, 84030-900, Brasil \\ bielinskilf@gmail.com}} 

\author[G. G. La Guardia]{\vspace{0.3cm} \\ Giuliano G. La Guardia  \\ {\tiny Universidade Estadual de Ponta Grossa\\ Ponta Grossa, PR, 84030-900, Brasil \\ gguardia@uepg.br}}

\author[J. Q. Chagas]{\vspace{0.3cm} \\ Jocemar Q. Chagas \\ {\tiny Universidade Estadual de Ponta Grossa \\ Ponta Grossa, PR, 84030-900, Brasil \\ jocemarchagas@uepg.br}}

\date{}

\begin{document}
	\maketitle

	\begin{abstract}
		This paper is concerned with the study of the fractional finite sums theory. We present the classes of functions for which it is possible to characterize the constant related to the derivative of fractional sums (denominated by {\it essence} of a function). Examples of the essences of functions are presented, including the Bernoulli numbers, which are obtained as essences of the functions $f(x)=x^a$, $a \in \mathbb{N}^*$. Based on the new results, we can expand the comprehension of the Euler-Maclaurin summation formula for a real number instead of only considering natural numbers. Moreover, we utilize the concept of the essence of a function to propose a new regularization method for divergent series.
	\end{abstract}
	
	\vspace{0.5cm} \noindent 
	
	%
	%------ Section 1 - Introduction -----
	%
	\section{Introduction}\label{Section1}
	
	The succinct and elegant notation for summation, $\displaystyle \sum_{i=m}^{n} f(i)$ for $m, n \in \mathbb{N}$, created by L. Euler in 1755 \cite{EulerE212}, is well-known in the literature. A few years ago, some number theory researchers also dealt with the understanding of the meaning of the symbols $\displaystyle \sum_{i=x}^{y} f(i)$, in case when $x, y$ are real or even complex numbers (see \cite{Muller2005,Muller2010,Muller2011}). Nowadays, such a type of sums are called {\it fractional finite sums}. The first known example of a fractional finite sum, due to L. Euler \cite{EulerE713}, is the identity
	\begin{equation}\label{Eq01}
		\displaystyle
		\sum_{k=1}^{-1/2} \frac{1}{k} = -2 \ln(2).
	\end{equation}
	S. Ramanujan also wrote some ideas concerning fractional sums \cite{Ramanujan2012,Berndt1985} (see, for example, Chapter VI of the second notebook). However, a structured approach for the theory of fractional sums appeared only after 2005 \cite{Muller2005,Muller2010,Muller2011}, due to D. Schleicher and M. Müller. More recently, I. Alabdulmohsin have published a book covering the theory for fractional finite sums \cite{Alabdulmohsin2018}, which includes an approach for fractional sums of functions that can oscillate the signal. The theory of fractional finite sums was also investigated in \cite{Machado2014,Galvao2021,Chagas2021} and applied in \cite{Machado2016,Uzun2019}. Techniques of umbral calculus were utilized in \cite{Qian2021} in order to deal with an open problem exhibited in \cite{Muller2010}.
	
	In this paper, we present some contributions for the theory of fractional finite sums. In Section~\ref{Section2}, a motivation to study fractional sums based on a simple example is shown. In Section~\ref{Section3}, we review basic concepts with respect to such a theory established by M. Müller and D. Schleicher. In Section~\ref{Section4}, we present some new results, including the answer to an open question raised in \cite{Muller2010}. We exhibited several examples and applications in Section~\ref{Section5}, including: (i) a derivation of the Euler-Maclaurin summation formula for fractional finite sums, which is valid for a large class of functions; and (ii) a new possibility to regularize divergent series, an important issue in physics. Finally, in Section~\ref{Section6}, the final remarks are drawn, as well as some open questions are raised.
	
	%
	% ----- Section 2  -----
	%
	\section{A motivation for the theory of fractional finite sums}\label{Section2}
	
	Some summations have closed expressions, i.e., they can be evaluated by using a fixed quantity of operations. For example, the sum of the first $n$ positive integers, $\displaystyle S_1(n)=\sum_{i=1}^ni$, can be evaluated by $S_1(n)=\frac{n(n+1)}{2}$. We can think in replacing $n \in \mathbb{N}$ by $x \in \mathbb{R}$ or even by $z \in \mathbb{C}$, and the expressions $S_1(x)$ or $S_1(z)$ still have a clear meaning. We can then write
	\begin{equation*}
		\displaystyle
		\frac{\pi(\pi+1)}{2} = S_1(\pi) := \sum_{i=1}^\pi i \,,
	\end{equation*}
	meaning ``the sum of the first $\pi$ integers''. Analogously, it is possible to propose ``continued versions'' for other sums having closed expressions, as well as $\displaystyle \sum_{i=1}^ni^2, \sum_{i=1}^ni^3$ and $\displaystyle \sum_{i=1}^nq^{i}$.
	Unfortunately, not all summations have closed expressions (or, even if there exists, it is very difficult to compute it). For these types of summations, a different approach is necessary.
	
	An example of finite sums that does not have a closed expression are the {\it harmonic numbers}, given by $\displaystyle H(n)=\sum_{i=1}^n\frac{1}{i}$.
	In this section, we present an introductory approach to obtain an expression which evaluates the ``continuous'' harmonic numbers from 1 to a real (or complex) summation boundary value $x$ (or $z$).
	Given a fixed positive integer $n$, for each integer $m > n$, by considering
	\begin{equation*}
		\begin{array}{rcllcl}
			A(m) & := & \displaystyle \sum_{i=1}^m \frac{1}{i} & = \displaystyle \frac{1}{1} + \frac{1}{2} + \dots + \frac{1}{n} & \displaystyle + \frac{1}{n+1} + \dots  +  \frac{1}{m} &   \\
			B(m) & := & \displaystyle \sum_{i=1}^m \frac{1}{i+n} & = & \displaystyle \ \ \ \frac{1}{n + 1} + \dots + \frac{1}{m} & \displaystyle + \dots + \frac{1}{m+n} \,,
		\end{array}
	\end{equation*}
	we obtain
	\begin{equation*}
		\displaystyle
		A(m)-B(m)= \sum_{i=1}^m \left(\frac{1}{i}-\frac{1}{i+n}\right) = H(n) - \sum_{i=1}^n {\frac{1}{m+i}}.
	\end{equation*}
	However, since
	\begin{equation}\label{DESIGUALDADE 1}
		\displaystyle 
		L(m) = \frac{n}{m+n} = \sum_{i=1}^n \frac{1}{m+n} \leq \sum_{i=1}^n \frac{1}{m+i} \leq \sum_{i=1}^n \frac{1}{m+1} = \frac{n}{m+1} = U(m),
	\end{equation}
	taking $m \to \infty$ in \eqref{DESIGUALDADE 1}, we have $\displaystyle\lim_{m\to\infty} L(m)=\lim_{m\to\infty} U(m) = 0$.
	From the squeeze theorem, it follows that $\displaystyle \lim_{m \to \infty} \Bigl(\sum_{i=1}^n \frac{1}{m+i}\Bigr) = 0$.
	Therefore, taking $m \to \infty$ in $\displaystyle \bigl(A(m)-B(m)\bigr)$, we conclude that for every positive integer $n$, the following equation holds:
	\begin{equation}\label{GERAL}
		\displaystyle
		H(n) = \sum_{i=1}^{\infty} \left( \frac{1}{i} - \frac{1}{i+n} \right).
	\end{equation}
	
Note that, for each index $i$, the function $\displaystyle f_i(n) = \frac{1}{i} - \frac{1}{i+n}$ in \eqref{GERAL}
has also meaning if we consider $n \in \mathbb{C} \setminus \{-i\}$. Hence, we propose the following formula:
\begin{equation}\label{EXPANSAO}
	\displaystyle
	\mathcal{H}(z) = \sum_{i=1}^{z} \frac{1}{i} := \sum_{i=1}^{\infty} \left( \frac{1}{i} - \frac{1}{i+z} \right) %= \psi(x+1)+\gamma,
\end{equation}
which is valid for all $z \in \mathbb{C} \setminus {\mathbb{Z}}_{-}$.

%\clearpage
\begin{figure}[h]
	\centering
	\includegraphics[scale=0.75]{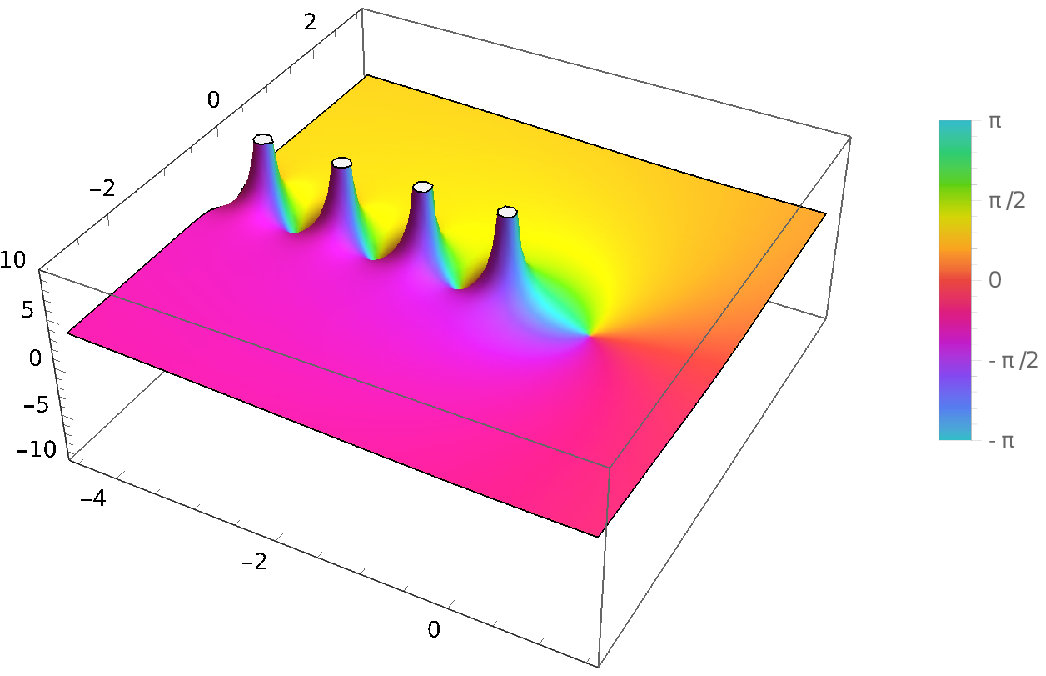}
	\caption{$\mathcal{H}(z)$ for $z\in\left[-4.5,1.5\right]\times \left[-3i,3i\right]$.}\label{Figure1}
\end{figure} 

In order to illustrate Formula~\eqref{EXPANSAO}, Figure~\ref{Figure1}
show $\mathcal{H}(z)$ for $z\in\left[-4.5,1.5\right]\times \left[-3i,3i\right]$.
For example, if $x=0$, we have $\displaystyle \mathcal{H}(0) = \sum_{i=1}^{\infty} \Bigl( \frac{1}{i} - \frac{1}{i+0} \Bigr) = 0$.
It is interesting to mention that the example given by L. Euler in \cite{EulerE713} was $\mathcal{H} \left( -\frac{1}{2} \right) = -2\ln(2)$. Formula~\eqref{EXPANSAO} is exactly the same that we obtain when evaluating the fundamental formula for the fractional sum established in \cite{Muller2011} (see Section~\ref{Section3}) for the function $\displaystyle f(\nu) = \frac{1}{\nu}$.
	
	%
	% ----- Section 3  -----
	%
	\section{Preliminaries}\label{Section3}
	
	In this section, we recall basic concepts of the theory of fractional finite sums established by M. Müller and D. Schleicher  \cite{Muller2005,Muller2010,Muller2011}.
	
	The axiomatic framework for the theory of fractional finite sums is composed by six axioms, presented by Müller and Schleicher in \cite{Muller2011}.
	For $x, y, z, s \in \mathbb{C}$ and $f, g$ complex functions (defined under some conditions which will be explained in the sequence), the following statements hold:
	
	\vspace{0.3cm}
	\noindent
	Axiom A1 (Continued summation):
	\begin{equation}\label{Axioma 1}
		\displaystyle
		\sum_{\nu=x}^{y} f(\nu) \ + \sum_{\nu=y+1}^{z} f(\nu) = \sum_{\nu = x}^{z} f(\nu).
	\end{equation}
	
	\noindent
	Axiom A2 (Translation invariance):
	\begin{equation}\label{Axioma 2}
		\displaystyle
		\sum_{\nu=x+s}^{y+s} f(\nu) = \sum_{\nu=x}^{y} f(\nu+s).
	\end{equation}
	
	\noindent
	Axiom A3 (Linearity):
	\begin{equation}\label{Axioma 3}
		\displaystyle
		\sum_{\nu=x}^{y} \bigl(\lambda f(\nu) + \mu \, g(\nu)\bigr) = \lambda \! \sum_{\nu=x}^{y} f(\nu) + \, \mu \! \sum_{\nu=x}^{y} g(\nu).
	\end{equation}
	
	\noindent
	Axiom A4 (Consistency with the classical sums):
	\begin{equation}\label{Axioma 4}
		\displaystyle
		\sum_{\nu=1}^{1} f(\nu) = f(1).
	\end{equation}
	
	\noindent
	Axiom A5 (Sums of monomials): for each $d \in \mathbb{N}$, the mapping
	\begin{equation}\label{Axioma 5}
		\displaystyle
		z \mapsto \sum_{\nu=1}^{z} \nu^d
	\end{equation}
	is holomorphic in $\mathbb{C}$.
	
	\vspace{0.3cm} \noindent
	Axiom A6 (Right shift continuity):
	if $\displaystyle \lim_{n \to \infty} f(z+n) = 0$ pointwise for any $z \in \mathbb{C}$, then
	\begin{equation}\label{Axioma 6}
		\displaystyle
		\lim_{n \to \infty} \sum_{\nu = x}^{y}  f(\nu+n) = 0.
	\end{equation}
	Moreover, if there exists a sequence ${(p_n)}_{n \in \mathbb{N}}$ of polynomials of fixed degree satisfying the condition
	$\left|f(z + n) - p_n (z + n) \right| \to 0$
	when $n \to \infty$, for all $z \in \mathbb{C}$, then
	\begin{equation}\label{Eq09}
		\displaystyle
		\bigg| \sum_{\nu=x}^{y} f(\nu + n) - \sum_{\nu=x}^{y} p_n (\nu + n) \biggr| \longrightarrow 0.
	\end{equation}
	
	\vspace{0.3cm}
	Axiom A6 has the following alternative version: if $\displaystyle \lim_{n \to \infty} f(z-n) = 0$ pointwise, for any $z \in \mathbb{C}$, then it follows that $\displaystyle \lim_{n \to \infty} \sum_{\nu=x}^{y}  f(\nu-n) = 0$. This alternative axiom is called ``left shift axiom'' or ``Axiom $\overleftarrow{\text{A6}}$''.
	
	By utilizing Axioms A1~-~A6, Müller and Schleicher developed a summation method valid for polynomials, as follows. Given a polynomial $p(z)$ with $z \in \mathbb{C}$, they considered $P : \mathbb{C} \rightarrow \mathbb{C}$ to be the unique polynomial satisfying the equation $P(z) - P(z - 1) = p(z)$ for all $z \in \mathbb{C}$ and the condition $P(0) = 0$. Thus, the unique definition for the fractional finite sums of polynomials $p(z)$, satisfying Axioms A1~-~A6, is given by~\cite{Muller2011}:
	\begin{equation}\label{Eq15}
		\displaystyle
		\sum_{\nu=x}^{y} p(\nu) =  P(y) - P(x-1) \,.
	\end{equation}
	Subsequently, they established, for real and complex functions, the (unique possible) definition for fractional finite sums, as follows~\cite{Muller2011}.
	
	\begin{definition}\label{fracsum} A function $f:U~\subset~\mathbb{C}~\to~\mathbb{C}$ is said to be fractional summable if the following conditions hold:
		\begin{itemize}
			\item[(i)] if $x \in U$, then $x + 1 \in U$;
			\item[(ii)] there exists a sequence ${(p_n)}_{n \in \mathbb{N}}$ of polynomials of degree less than or equal to $m$,  such that, for all $x \in U$, we have $\displaystyle \bigl|f(n+x) - p_n(n+x)\bigr| \to 0$ when $n \to + \infty$ (by convention, the null polynomial has degree $- \infty$); and
			\item[(iii)] for all $x, y \in U$, there exists the limit
			\begin{equation}\label{limit}
				\displaystyle
				\lim_{n \to \infty} \bigg( \sum_{\nu=n+x}^{n+y} p_n(\nu) + \sum_{\nu=1}^{n} \bigl(f(\nu+x-1) - f(\nu+y)\bigr) \biggr),
			\end{equation} where the fractional sum $\displaystyle \sum p_n$ is in the sense of Eq.~\eqref{Eq15}.
		\end{itemize}
	\end{definition}
	
	For fractional summable functions $f:U\subset\mathbb{C} \to \mathbb{C}$ such that $f(n+z)\to 0$ as $n\to\infty$ for all $z\in U$, M\"uller and Schleicher~\cite{Muller2011} shown that:
	\begin{equation}\label{Eq33}
		\displaystyle
		\sum_{\nu=x}^{y} f(\nu) =  \sum_{\nu = 1}^{\infty} \bigl( f(\nu+x-1) - f(\nu+y) \bigr) \, .
	\end{equation}
	
	They utilized the notation $\displaystyle {\sum_x^y \hspace{-0.45cm}\rightarrow f}$ to represent the limit in \eqref{limit}, and called it by fractional sum of $f$~\cite{Muller2011}. 
	In this article, these symbols have a slightly different use (see Definition~\ref{Definition4.04}).
	
	I. Alabdulmohsin developed the same theory (see \cite{Alabdulmohsin2018}) by means of only two axioms, namely, Axiom A1 and Axiom A4, rewritten as
	\begin{equation}\label{Eq06alt}
		\displaystyle
		\sum_{\nu=z}^{z} f(\nu) = f(z)
	\end{equation}
	for all $z \in U \subset \mathbb{C}$, where $f$ is fractional summable.	
	
	%
	% ----- Section 4  -----
	%
	\section{New results}\label{Section4}
	
	We present here some new contributions on the theory of fractional finite sums. One of our purposes is to answer an open problem raised in \cite{Muller2010}. More precisely, our main purposes are to establish conditions on a function $f$ in order to satisfy the relation
	\begin{equation}\label{Eq4.01}
		\displaystyle
		\frac{d}{dx} \sum_{\nu=1}^{x}\hspace{-0.525cm}\rightarrow f(\nu) \overset{?}{=} c_f
		+ \sum_{\nu=1}^{x}\hspace{-0.525cm}\rightarrow f'(\nu) \, ,
	\end{equation}
	as well as to define the non-arbitrary constant $c_f$, providing methods for evaluating it in some cases.
	
	\subsection{On Axioms, sets, and immediate consequences}
	
	We consider the complete axiomatic framework proposed by M. Müller and D. Schleicher \cite{Muller2011} (see Section \ref{Section3}), where Axiom A5 is replaced by Axiom A5* (given in the sequence) and A6 consists only in the right (but not in the left) shift continuity. In Proposition~\ref{Proposition4.01}, we present some immediate and useful consequences of Axioms A1, A2 and A4. We also introduce the concept of a summable set and define some other sets of functions, suitable for the development of the next subsections.
	
	\noindent{Axiom A5* (Sums of holomorphic functions) If $f:\mathbb{C}\to\mathbb{C}$ is holomorphic, then the function}
	\begin{equation*}
		\displaystyle
		x \mapsto \sum_{k=1}^{x}\hspace{-0.5cm}\rightarrow f(k)
	\end{equation*}
	is also holomorphic.
	
	Since every polynomial is holomorphic, Axiom A5* clearly implies Axiom A5.
	
	\begin{definition}\label{Definition4.02}
		For a set $A$, we denote by $\wp(A)$ the set of all subsets of $A$. Let $\mathbb{C}^\mathbb{C}$ the set of all functions $f:\mathbb{C}\to\mathbb{C}$. Define the functions $\mathbf{L},\mathbf{T},\mathbf{S}:\wp\left(\mathbb{C}^\mathbb{C}\right)\to\wp\left(\mathbb{C}^\mathbb{C}\right)$, respectively, by
		\begin{equation*}
			\begin{array}{rcl}
				\mathbf{L}(X) & := & \left\{f \in \mathbb{C}^\mathbb{C} \, | \, f \text{ is a linear combination of functions in } X \right\}, \\
				\mathbf{T}(X) & := & \left\{f \in \mathbb{C}^\mathbb{C} \, | \, \exists \ (g,c) \in X\times\mathbb{C}, \, \forall \ \, x \in \mathbb{C}, \, f(x)=g(x+c) \right\}, \ \text{and} \\
				\mathbf{S} & := & \mathbf{L}\circ\mathbf{T}.
			\end{array}
		\end{equation*}
	\end{definition}
	
	\begin{definition}\label{Definition4.05}
		Let $\mathcal{F}$ be the set of the functions $f:\mathbb{C}\to\mathbb{C}$ that can be expressed as a linear combination of holomorphic and fractional summable functions.
	\end{definition}
	
	\begin{definition}\label{Definition4.04}
		If $\mathcal{U}\subset\mathcal{F}$ is such that there exists a unique application
		\begin{equation}\label{Eq4.09}
			\begin{split}
				\sum\hspace{-0.45cm}\rightarrow\hspace{0.1cm} : \mathbf{S}(\mathcal{U})\times\mathbb{C}^2 & \to \mathbb{C} \\
				(f,x,y) &\mapsto \sum\hspace{-0.45cm}\rightarrow(f,x,y) =  \sum_{k=x}^y\hspace{-0.5cm}\rightarrow f(k)
			\end{split}
		\end{equation}
		satisfying axioms
		A1$_{\mathcal{U}}$, A2$_\mathcal{U}$, A3$_\mathcal{U}$, A4$_\mathcal{U}$, A5$^*_\mathcal{U}$ and A6$_\mathcal{U}$, then $\mathcal{U}$ is called a {\it summable set} and $\displaystyle\sum\hspace{-0.45cm}\rightarrow\hspace{0.1cm}$ is the {\it fractional sum relative to $\mathcal{U}$}.
		
		\begin{remark}
			Note that a function can belong to a summable set even though it is not fractional summable; see for instance, Example~\ref {Example5.05} for $z=\pi i$.
		\end{remark}
		
		More generally, if we have a scheme of axioms depending on a set $X\subset \mathbb{C}^{\mathbb{C}}$, that is, an scheme $\mathcal{A}$ given by $\mathcal{A}(X)=(\mathcal{A}_1(X), \mathcal{A}_2(X),...,\mathcal{A}_n(X))$, and $\mathcal{U}\subset\mathbb{C}^\mathbb{C}$ is such that there exists a unique application $\displaystyle\sum\hspace{-0.45cm}\rightarrow\hspace{0.1cm} : \mathbf{S}(\mathcal{U})\times\mathbb{C}^2 \to \mathbb{C}$ satisfying such axioms, then $\mathcal{U}$ is called a {\it summable set in scheme of axioms $\mathcal{A}$}.
		Finally, we define $\mathfrak{S}:=\left\{X\subset\mathcal{F}\,|\,X\text{ is a summable set}\right\}$.
	\end{definition}
	
	\begin{remark}\label{Remark04}
		For each set $\mathcal{U}$, the set of functions in question will be denoted by $\mathbf{S}(\mathcal{U})$. That is, each set $\mathcal{U}$ generates their respective axioms A1$_\mathcal{U}$, A2$_\mathcal{U}$, A3$_\mathcal{U}$, etc. For example, Axiom A2$_\mathcal{U}$ says that $\forall \ (f,s,x,y) \in \mathbf{S}(\mathcal{U})\times\mathbb{C}^3, \, \displaystyle \sum\hspace{-0.45cm}\rightarrow (f,x+s,y+s) = \sum\hspace{-0.45cm}\rightarrow (h,x,y)$, where $\forall \ z \in \mathbb{C}, \, h(z):=f(z+s)$. Axiom A3$_\mathcal{U}$ says that $\forall \ (f,g,\lambda,\mu,x,y) \in \mathbf{S}(\mathcal{U})^2\times\mathbb{C}^4, \, \displaystyle \sum\hspace{-0.45cm}\rightarrow (\lambda f+\mu g,x,y) = \lambda\sum\hspace{-0.45cm}\rightarrow (f,x,y) + \mu \sum\hspace{-0.45cm}\rightarrow (g,x,y)$. Axiom A5$_{\mathcal{U}}$ says that if $f\in\mathcal{U}$ is holomorphic, then it is holomorphic the mapping $\displaystyle x \mapsto \sum\hspace{-0.45cm}\rightarrow (f,1,x)$. Note that $h\in\mathbf{T}(\mathcal{U})\subset\mathbf{S}(\mathcal{U})$ and $\lambda f + \mu g\in\mathbf{L}(\mathcal{U})\subset\mathbf{S}(\mathcal{U})$.
	\end{remark}
	
	\begin{remark}\label{globalmeaning}
		Note that, a piori, the symbol $\displaystyle\sum\hspace{-0.45cm}\rightarrow\hspace{0.1cm}$ has no ``global meaning'': let $\mathcal{U},\mathcal{V}$ two summable sets, let $^\mathcal{U}\hspace{-0.05cm}\displaystyle\sum\hspace{-0.45cm}\rightarrow$ and $^\mathcal{V}\hspace{-0.05cm}\displaystyle\sum\hspace{-0.45cm}\rightarrow$ the fractional sums relative to the sets $\mathcal{U}$ and $\mathcal{V}$, and let $f \in \mathbf{S}\left(\mathcal{U} \cap \mathcal{V}\right)$. Then, we cannot guarantee that $\displaystyle^\mathcal{U}\hspace{-0.05cm}\sum\hspace{-0.45cm}\rightarrow\hspace{0.1cm}f = ^\mathcal{V}\hspace{-0.05cm}\displaystyle\sum\hspace{-0.45cm}\rightarrow\hspace{0.1cm}f$ (see open question (ii)).
	\end{remark}
	
	\begin{remark}\label{omit}
		When it is clear from the context, as in Propositions~\ref{notsummable1} and~\ref{Proposition4.01}, we omit the set in which the fractional sum is related.
	\end{remark}
	
	Müller and Schleischer \cite{Muller2011} show us that the set $\mathbb{C}[x]$ of the complex polynomials is a summable set; therefore, $\mathfrak{S}\neq\emptyset$.
	One of the aims of this paper consists of the searching for the summable sets and to compute the fractional sums of their elements.
	Another examples of summable sets are:
	\begin{itemize}
		%\item[i) ] The set of all the functions $f:\mathbb{C}\to\mathbb{C}$ that are fractional summable (obviously including the polynomials);
		\item[i) ] $\{x\mapsto u^x\,|\,u\in\mathbb{C}\setminus\{0,1\}\}$;
		\item[ii) ] $\mathbb{C}[x]\cup\{x\mapsto u^x\,|\,u\in\mathbb{C}\setminus\{0,1\}\}$;
		\item[iii) ] $\{P\in\mathbb{C}[x]\,|\,P\text{ has degree less than or equal to } c\}$, where $c$ is a positive integer;
		\item[iv) ]$\left\{f\in\mathbb{C}^\mathbb{C}\,|\,f(0)=0,\exists\,n\in\mathbb{N}^*,\forall \,x\in\mathbb{C}^*,f(x)=\frac{1}{x^n}\right\}$.
	\end{itemize}
	
	Unfortunately, not all holomorphic functions can be ``fractionally summed'', as shown in Proposition~\ref{notsummable1}. In fact, such a function cannot be summable in any scheme of axioms that includes axioms A1-A4. This is an example of a ``extremely well-comported'' holomorphic function that cannot be ``fractional-summed''.
	
	\begin{proposition}\label{notsummable1}
		For every odd integer $k$, any set containing the function $x\mapsto e^{2\pi i k x}$ is not a summable set. More generally, if $f:\mathbb{C} \to \mathbb{C}$ is such that $f \left( \frac{1}{2} \right) \neq 0$ and for all $x \in \mathbb{C}$ satisfies $f \left (x + \frac{1}{2} \right) = -f(x)$, then any set containing $f$ is not a summable set in every scheme of axioms including Axioms A1, A2, A3 and A4.
	\end{proposition}
	\begin{proof}
		Suppose that $\mathcal{U}$ is a  summable set and let $c=\displaystyle\sum_{k=1}^{1/2}\hspace{-0.5cm}\rightarrow \hspace{0.05cm}f(k)$. Then
		\begin{equation*}
			\begin{array}{rrclr}
				& c+\displaystyle\sum_{k=3/2}^{1}\hspace{-0.65cm}\rightarrow\hspace{0.15cm} f(k) & = &\displaystyle\sum_{k=1}^{1/2}\hspace{-0.5cm}\rightarrow\hspace{0.05cm} f(k)+\sum_{k=3/2}^{1}\hspace{-0.65cm}\rightarrow\hspace{0.15cm} f(k) & 
				\\
				\Rightarrow & \quad c+\displaystyle\sum_{k=1}^{1/2}\hspace{-0.5cm}\rightarrow\hspace{0.05cm} f\left(k+\frac{1}{2}\right)&=&\displaystyle\sum_{k=1}^{1}\hspace{-0.5cm}\rightarrow\hspace{0.05cm} f(k)&\qquad\text{(from A1 and A2)}
				\\
				\Rightarrow & c-\displaystyle\sum_{k=1}^{1/2}\hspace{-0.5cm}\rightarrow\hspace{0.05cm} f\left(k\right)&=&f(1)&\qquad\text{(from A3 and A4)}
				\\
				\Rightarrow &  f(1) & = & 0. &
			\end{array}
		\end{equation*}
		But $f\left(\frac{1}{2}+\frac{1}{2}\right)=-f\left(\frac{1}{2}\right)\Rightarrow f(1)=-f\left(\frac{1}{2}\right)\Rightarrow f\left(\frac{1}{2}\right) = 0$, a contradiction. For the particular case, just note that $e^{2\pi i k \frac{1}{2}}\neq 0$ and $e^{2\pi i k \left(x+\frac{1}{2}\right)}=e^{2\pi i k x}e^{2\pi i k \frac{1}{2}}=-e^{2\pi i k x}$.
	\end{proof}
	
	In Proposition~\ref{Proposition4.01}, we present some immediate and useful consequences derived from Axioms A1-A4.
	
	\begin{proposition}\label{Proposition4.01}
		Let $\mathcal{U}$ a summable set, $f\in \mathbf{S}(\mathcal{U})$ and $x,y \in \mathbb{C}$. Then, the following hold:
		\begin{itemize}
			\item[i) ] Consistency with the classical sums:
			\begin{equation}\label{Eq4.02}
				\displaystyle
				\sum_{k=x}^{x}\hspace{-0.5cm}\rightarrow f(k) = \ f(x) .
			\end{equation}
			
			\item[ii) ] Continued summation from 1:
			\begin{equation}\label{1Convention}
				\displaystyle
				\sum_{k=x}^{y}\hspace{-0.5cm}\rightarrow f(k) = \sum_{k=1}^{y}\hspace{-0.5cm}\rightarrow f(k) - \sum_{k=1}^{x-1}\hspace{-0.5cm}\rightarrow f(k) .
			\end{equation}
			\noindent
			\item[iii) ] Generalized empty sum:
			\begin{equation}\label{EmptySum}
				\displaystyle
				\sum_{k=x}^{x-1}\hspace{-0.5cm}\rightarrow f(k) = 0 .
			\end{equation}
			
			\item[iv) ] Interpolation of classical sums
			\begin{equation}\label{Interpolation}
				\displaystyle
				\sum_{k=x}^{y}\hspace{-0.5cm}\rightarrow f(k) = \sum_{k=a}^{y-x+a}f(k+x-a), \quad \text{if }a,y-x\in\mathbb{N} .
			\end{equation}
			
			\item[v) ] Opposite sum
			\begin{equation}\label{Opposite sum}
				\displaystyle
				\sum_{k=1}^{-x}\hspace{-0.5cm}\rightarrow f(k) = -\sum_{k=1}^{x}\hspace{-0.5cm}\rightarrow f(k-x) .
			\end{equation}
		\end{itemize}
	\end{proposition}
	\begin{proof}
		For Items (i) and (ii), apply A2, A4 and A1, respectively. Item (iii) follows from Item (ii) by taking $y=x-1$. For Item (iv), just apply A2, A1, and then Item (i). Finally, for Item (v), apply A2 and Item (ii).
	\end{proof}
	
	Since it is possible to replace any fractional summation by equivalent sums beginning at 1, in the following we consider only sums beginning at $k=1$. When it is convenient, we write only $\displaystyle\sum^x\hspace{-0.45cm}\rightarrow f$ instead of writing $\displaystyle\sum_{k=1}^x\hspace{-0.5cm}\rightarrow f(k)$. Moreover, the idea of the empty sums, given in Proposition~\ref{Proposition4.01}, already appears in \cite{Alabdulmohsin2018}. In this work is fundamental the concept of empty sum.
	
	\begin{remark}\label{Remark02}
		It is interesting to note the follow. Suppose that our scheme of axioms is (A1, A2, A3, A4, A5, A6). For $\mathcal{U}=\{x\mapsto u^x\}$ with $|u|<1$ and $u\neq 0$, we could simply apply Axiom A6 to evaluate the fractional sums over all the complex plane directly. However, the set $\mathcal{V}=\{x\mapsto u^x\}$, where $|u|\geq1$ and $u\neq 1$, would not be a summable set. In this case, we have that:
		\begin{equation*}
			\begin{array}{rcl}
				\mathbf{S}(\mathcal{V}) & = &\displaystyle\left\{\left.x\mapsto \sum_{m=1}^n\alpha_m u^{x+\beta_m}\,\right|\,n\in\mathbb{N}^*,(\alpha,\beta)\in\mathbb{C}^{2n}\right\}
				\\
				& = & \displaystyle\left\{\left.x\mapsto \sum_{m=1}^n\omega_m u^{x}\,\right|\,n\in\mathbb{N}^*,\omega\in\mathbb{C}^{n}\right\}
				\\
				& = & \displaystyle\left\{\left.x\mapsto a u^{x}\,\right|\,a\in\mathbb{C}\right\}.
			\end{array}
		\end{equation*}
		In fact, for $p,q\in\mathbb{N}^*$, we have
		\begin{equation*}
			\begin{array}{rrclr}
				& \displaystyle\sum_{l=1}^q\sum_{k=1+(l-1)\frac{p}{q}}^{l\frac{p}{q}}\hspace{-1.0cm}\rightarrow\hspace{0.5cm}u^k & \overset{A1}{=} &\displaystyle \sum_{k=1}^{p}\hspace{-0.5cm}\rightarrow u^k&
				\\
				\Rightarrow & \quad \displaystyle\sum_{l=1}^q\sum_{k=1}^{\frac{p}{q}}\hspace{-0.5cm}\rightarrow u^{k+(l-1)\frac{p}{q}} & {=} &\displaystyle \sum_{k=1}^{p} u^k&\quad\text{Axiom A2 and Eq.~\eqref{Interpolation}}
				\\
				\Rightarrow & \displaystyle\sum_{l=1}^q u^{(l-1)\frac{p}{q}}\displaystyle\sum_{k=1}^{\frac{p}{q}}\hspace{-0.5cm}\rightarrow u^k & {=} &\displaystyle \dfrac{u\left(u^p-1\right)}{u-1}&\quad\text{Axiom A3}
				\\
				\Rightarrow & \displaystyle\sum_{k=1}^{\frac{p}{q}}\hspace{-0.5cm}\rightarrow u^k& = &\dfrac{u\left(\cancel{u^p-1}\right)}{u-1}\cdot\dfrac{1}{\left(\frac{\cancel{1\left(\left(u^\frac{p}{q}\right)^q-1\right)}}{\left(u^{\frac{p}{q}}\right)-1}\right)}&
				\\
				\Rightarrow & \displaystyle\sum_{k=1}^{\frac{p}{q}}\hspace{-0.5cm}\rightarrow u^k& = & \dfrac{u\left(u^\frac{p}{q}-1\right)}{u-1}
				\\
				\Rightarrow & \displaystyle\sum_{k=1}^{\frac{p}{q}}\hspace{-0.5cm}\rightarrow au^k& = & a\dfrac{u\left(u^\frac{p}{q}-1\right)}{u-1}.
			\end{array}
		\end{equation*}
		
		In case of negative rationals, note that
		\begin{equation*}
			\begin{array}{rclr}
				\displaystyle\sum_{k=1}^{-\frac{p}{q}}\hspace{-0.5cm}\rightarrow au^k & = & - \displaystyle\sum_{k=1}^{\frac{p}{q}}\hspace{-0.5cm}\rightarrow au^{k-\frac{p}{q}} & \quad \text{by \eqref{Opposite sum}}
				\\
				& = &\displaystyle -u^{-\frac{p}{q}}\sum_{k=1}^{\frac{p}{q}}\hspace{-0.5cm}\rightarrow au^{k} & \quad \text{Axiom A3}
				\\
				& = & -u^{-\frac{p}{q}}\left(a \dfrac{u\left(u^\frac{p}{q}-1\right)}{u-1}\right)
				\\
				& = &
				a \dfrac{u\left(u^{-\frac{p}{q}}-1\right)}{u-1}.&
			\end{array}
		\end{equation*}
		Finally, for $0$, recall $\eqref{EmptySum}$ and note that $0 = a \dfrac{u\left(u^0-1\right)}{u-1}$. Therefore,
		\begin{equation}\label{expsum}
			\displaystyle
			\sum_{k=1}^{r}\hspace{-0.5cm}\rightarrow au^k = a\dfrac{u\left(u^r-1\right)}{u-1}, \quad \forall \ r \in \mathbb{Q}.
		\end{equation}
		
		However, we do not have perspective to evaluate this fractional sum for $r\notin\mathbb{Q}$, and only Axiom A5$^*$ permits us to evaluate the missing points. Note that the left shift axiom only solves these problems for $|u|>1$, but not for $|u|=1$; hence, one cannot evaluate, for example, $\displaystyle\sum_{k=1}^x\hspace{-0.5cm}\rightarrow i^k:=\sum_{k=1}^x\hspace{-0.5cm}\rightarrow e^{(\pi i/2)k}$.
	\end{remark}
	
	\begin{remark}\label{Remark01} 
		The above considerations are justifications for not considering $f\notin\mathcal{F}$; apparently, Axioms A1, A2, A3 and A4, in the best scenario, can be used to determine $\displaystyle\sum_x^y\hspace{-0.45cm}\rightarrow f=\displaystyle\sum^y\hspace{-0.45cm}\rightarrow f-\displaystyle\sum^{x-1}\hspace{-0.5cm}\rightarrow f$ only for $y-x\in\mathbb{Q}$, but not for $y-x\notin\mathbb{Q}$, and we cannot apply Axioms A5$^*$ and A6 in any way to solve the missing case.
		
		Because of this fact, we consider two distinct functions $S,P:\mathbb{C}\times\mathbb{C}\to\mathbb{C}$ such that
		\begin{equation*}
			\Bigl(\sum\hspace{-0.45cm}\rightarrow f \Bigr)\Big|_{Q}=S\big|_{Q}=P\big|_{Q}
		\end{equation*}
		only for $Q\subset{\{(x,y)\in\mathbb{C}^2\,|\,y-x\in\mathbb{Q}\}}$, but not for $Q=\mathbb{C}$, which break our uniqueness assumption.
	\end{remark}
	
	\subsection{Essence of a function}
	
	The first step to solve the open question left in \cite{Muller2011}, about conditions under which the relation \eqref{Eq4.01} holds, is to characterize uniquely a constant relative to each summable function $f$. With this purpose in mind, we define such a constant as follows.
	
	\begin{definition}[Essence of a function]\label{Definition4.03}
		Let $\mathcal{U}$ be a summable set, $\displaystyle\sum\hspace{-0.45cm}\rightarrow$ be the fractional sum relative to $\mathcal{U}$, and $f\in\mathbf{S}(\mathcal{U})$. The {\it essence of $f$ relative to $\mathcal{U}$} is a constant given by
		\begin{equation}\label{Eq4.13}
			\mathrm{ess}(f):=\lim_{h\to0} \frac{1}{h} \sum^h\hspace{-0.45cm}\rightarrow f,
		\end{equation}
		if the limit exists. In this case we say that $f$ is a {\it function with essence in $\mathcal{U}$}.
	\end{definition}
	
	\begin{remark}
		Remarks~\ref{globalmeaning} and~\ref{omit} should be considered here in analogous way.
	\end{remark}
	
	It is easy to verify that the essence is linear in $\mathbf{S}(\mathcal{U})$.
	
	Given a function $f$, we will denote $\mathrm{ess}(f)$ by $\mathrm{ess}(x\mapsto f(x))$. For example, $\mathrm{ess}(y\mapsto y^2-2)$ means the same of writting $\mathrm{ess}(g)$ for $g(y)=y^2-2$.
	
	\begin{proposition}[Characterization of the essence]\label{Proposition4.02}
		If $f$ is a function with essence, we have
		\begin{align}\label{Eq4.14}
			\displaystyle
			\mathrm{ess}(f) = \left. \frac{d}{dx} \Bigl( \sum^{x}\hspace{-0.45cm}\rightarrow f \Bigr) \right|_{x=0} .
		\end{align}
	\end{proposition}
	
	\begin{proof}
		It is sufficient to apply the derivative with respect to $x$ at $\displaystyle \sum^x\hspace{-0.45cm}\rightarrow f$, and evaluating it at $x=0$:
		\begin{align*}%\label{Eq4.15}
			\displaystyle
			\left. \frac{d}{dx} \Bigl( \sum^x\hspace{-0.45cm}\rightarrow f \Bigr) \right|_{x=0}
			& =
			\displaystyle
			\left. \lim_{h\to0} \, \frac{1}{h} \left( \sum_{k=1}^{x+h}\hspace{-0.5cm}\rightarrow f(k) - \sum_{k=1}^x\hspace{-0.5cm}\rightarrow f(k) \right) \right|_{x=0}  \\
			& =
			\displaystyle
			\lim_{h\to0} \, \frac{1}{h} \left( \sum_{k=1}^{0+h}\hspace{-0.5cm}\rightarrow f(k) - \sum_{k=1}^0\hspace{-0.5cm}\rightarrow f(k) \right) \\
			& =
			\displaystyle \lim_{h\to0} \, \frac{1}{h} \sum_{k=1}^{h}\hspace{-0.5cm}\rightarrow f(k) \\
			& =
			\displaystyle
			\mathrm{ess}(f)\,.
		\end{align*}
	\end{proof}
	
	Axiom A5$^*$ implies that every holomorphic function that belongs to a summable set also has essence in that set.
	
	\begin{example}\label{Example4.02}
		In Table~\ref{tab1}, we display some examples of functions and their respective essences.
	\end{example} 
	
	\begin{remark}
		All functions in Table~\ref{tab1} have fractional sum in a summable set $\mathcal{U}$ that can be calculated by only requiring that $f \in \mathbf{S}(\mathcal{U})$ (that is, we just use the functions in $\mathbf{S}(\{f\})$ to do the calculations), so that these fractional sums (and, consequently, their respective essences) are unique, independently of the summable set they belong to.
	\end{remark}
	
	% ----- Table 1 -----
	\begin{table}[!h]\label{tab1}
		\center
		\small{{\bf Table~\ref{tab1}} - Essences for some selected functions.} \\
		\begin{tabular}{clclc}
			\hline
			\qquad & $f$ & \qquad & $\mathrm{ess}(f)$ & \qquad \\ \hline
			\qquad & $x \mapsto z$ & \ & $z$ & \qquad \\
			\qquad & $x \mapsto x^a,\,a\in\mathbb{Z}^*$ & \ & $-a\zeta(1-a)$ & \qquad \\
			%x \mapsto x^{w},\,w\neq 0 & \ & -w\cdot\zeta(1-w)
			%\\
			\qquad & $x \mapsto e^{zx}, \ e^z\neq 1$ & \ & $\dfrac{ze^z}{e^z-1}$ & \qquad \\ 
			\qquad & $x \mapsto e^{z x}x, \ e^z\neq 1$ & \ & $\dfrac{e^z}{e^z-1}\left(1-\dfrac{z}{e^z-1}\right)$ & \qquad \\ 
			\qquad & $x \mapsto \ln x$ & \ & $-\gamma$ & \qquad \\
			\qquad & $x\mapsto \dfrac{1}{x(x+1)}$ & \ & $1$ & \qquad \\
			%\\  x \mapsto e^{\pi i x}x & \ & \frac{\pi i}{4} + \frac{1}{2} \\ 
			\hline
		\end{tabular}
	\end{table} 	
	
	\begin{remark}
		Note that  $\gamma$ is the Euler-Mascheroni constant.
	\end{remark}
	
	\begin{remark}\label{somenotation}
		Some of the above functions are not defined throughout $\mathbb{C}$. In this article, every time we refer to a function $f$ that can be defined in $\mathbb{C}-X$, we will actually consider the extension of $f$ by zero. That is, the function $\overline{f}:\mathbb{C}\to\mathbb{C}$ such that $\overline{f}|_{\mathbb{C}-X}=f$ and $\overline{f}(x)=0$ for all $x\in X$.
		For example, when we say $f(x)=\dfrac{1}{x}$ or $f(x)=x^{-1}$, we are referring to the function $f:\mathbb{C}\to\mathbb{C}$ such that $f|_{\mathbb{C}-\{0\}}(x)=\dfrac{1}{x}$ and $f(0)=0$.
	\end{remark}
	
	\begin{theorem}\label{Theorem4.01}
		Let $a\in\mathbb{N}$. Then, in any summable set $\mathcal{U}$ such that $\mathbf{S}(\mathcal{U})$ contains the function $x\mapsto x^a$, it follows that $\mathrm{ess}\left(x\mapsto x^a\right)=B_a$, where $B_a$ are the Bernoulli numbers with $B_1=\frac{1}{2}$.
	\end{theorem}
	\begin{proof}
		For $a = 0$, $x\mapsto x^a$ is the constant function $x\mapsto 1$. Then, $\mathrm{ess}(x\mapsto x^a)=\displaystyle\lim_{h\to0}\dfrac{1}{h}\sum^h\hspace{-0.45cm}\rightarrow1=\lim_{h\to0}\dfrac{1}{h}h=1=B_0$.
		
		For $a>0$, (see \cite{Bernoulli1713,Arakawa2014}) for all $(n,a)\in\mathbb{N}^*\times\mathbb{N}^*$, we have
		\begin{equation*}%\label{Eq4.34}
			\displaystyle
			\sum_{k=1}^n k^a = n^a + \sum_{k=0}^a \binom{a}{k} B_{a-k}^{-} \dfrac{n^{k+1}}{k+1},
		\end{equation*}
		where $B_1^-=-\frac{1}{2}$ and $B_i^-=B_i$ for $i\neq 1$. From  Item (iv) of Proposition~\ref{Proposition4.01}, it follows that for all $(x,a)\in\mathbb{N}^*\times\mathbb{N}^*$, we have
		\begin{equation}\label{Eq4.50}
			\displaystyle
			\sum_{k=1}^{x}\hspace{-0.5cm}\rightarrow k^a = x^a + \sum_{k=0}^{a} \binom{a}{k} B_{a-k}^{-} \dfrac{x^{k+1}}{k+1}.
		\end{equation}
		
		The second member of Eq.~\eqref{Eq4.50} is a polynomial in $x$ of degree $a+1$, and in \cite{Muller2011}, the same is true for the first member. But if $P$ and $Q$ are polynomials of degree $m$ and $P(x)=Q(x)$ for $m+1$ distinct values of $x$, then $P = Q$, so the Eq.~\eqref{Eq4.50} holds for all $x\in\mathbb{C}$. Therefore, assuming that $a\geq 2$, we can differentiate both sides of Eq.~\eqref{Eq4.50} at $0$:
		\begin{equation*}
			\begin{array}{rcl}
				\displaystyle\mathrm{ess}(x\mapsto x^a) &=& \displaystyle\left.ax^{a-1}+ \sum_{k=0}^a\binom{a}{k}B_{a-k}^-x^k\right|_{x=0}
				\\
				&=& 0 + \displaystyle\binom{a}{0}B_{a-0}^-\cdot 1
				\\
				&=& B_a.
			\end{array}
		\end{equation*}
		The case $a=1$ can be easily verified. Therefore, $\mathrm{ess}(x\mapsto x^a)=B_a=-a\zeta(1-a)$ for all $a\in\mathbb{N}^*$ \cite[p. 807]{Abramowitz1964}. %Eqs.~(23.2.11),(23.2.14),(23.2.15)
		Moreover, by linearity, if $\mathbb{C}[x]\subset\mathbf{S}(\mathcal{U})$ and $P(x) = c_0 + \displaystyle\sum_{k=1}^nc_kx^k$, then $\displaystyle\mathrm{ess}(P)=\sum_{k=0}^nc_kB_k$.
	\end{proof}
	
	Proposition~\ref{Proposition4.02} presents a characterization of the essence of a function $f$; however, in much cases, such a characterization is not suitable for evaluating the essence, since the function $\displaystyle x\mapsto \sum^{x}\hspace{-0.45cm}\rightarrow f$ cannot be known a priori. In Section~\ref{Section6}, we will exhibit methods that can be utilized to compute the essence of some functions.
	
	\subsection{{The main results}}
	
	We present here our main results, in which we determine some classes of functions for which it is possible to obtain its fractional summation, and the result obtained agree with Eq.~\eqref{Eq4.01}, where the constant $c_f$ is replaced by the essence of the function. More specifically, in Theorem~\ref{Theorem4.02}, we show that Eq.~\eqref{Eq4.19} is true for holomorphic functions belonging to a summable set satisfying some hypothesis which will be detailed in the following. In Theorem~\ref{tayloress} we expand the result to get a new one. In Theorem~\ref{Theorem4.04}, we show that the result is also valid for some functions which are not covered by Theorems~\ref{Theorem4.02}-\ref{tayloress}. Finally, in Theorem~\ref{Theorem5}, we present some identities about the essence of the translation of a function by a complex number $y$.
	
	\begin{theorem}\label{Theorem4.02}
		Let $f:\mathbb{C}\to\mathbb{C}$ an holomorphic function. If there exist $2n+1$ functions $u,v_1,\dots,v_n,w_1,\dots,w_n:\mathbb{C}\to\mathbb{C}$ satisfying
		\begin{itemize}
			\item[i) ]$\forall \ m\in\{1,\dots,n\}, \ \displaystyle \lim_{h\to 0}\left(\dfrac{v_m(h)}{h}\right)=0$,
			\item[ii) ]$\displaystyle \forall \ (x,h) \in \mathbb{C}\times\mathbb{C}, \ f(x+h)-f(x)=u(h)f'(x)+\sum_{m=1}^nv_m(h)w_m(x)$,
		\end{itemize}
		and if there exists a summable set $\mathcal{U}$ such that $\{f, f', w_1,\dots,w_n\}\subset \mathbf{S}(\mathcal{U})$,
		then it follows that
		\begin{equation}\label{Eq4.19}
			\displaystyle
			\frac{d}{dx} \left( \sum_{k=1}^{x}\hspace{-0.5cm}\rightarrow f(k) \right) = \mathrm{ess}(f) + \sum_{k=1}^{x}\hspace{-0.525cm}\rightarrow f'(k).
		\end{equation}
		
	\end{theorem}
	
	%\begin{remark}
	%Note that Eq.~\eqref{Eq4.19} is the same as Eq.~\eqref{Eq4.01} with $c_f$ replaced by $\mathrm{ess}(f)$.
	%\end{remark}
	
	\begin{proof}
		Let $f$ satisfying the hypotheses for Theorem~\ref{Theorem4.02}. So, for all $(x,h)\in\mathbb{C}\times\left(\mathbb{C}\setminus\{0\}\right)$, we have
		\begin{equation*}%\label{Eq4.22}
			\begin{array}{rcl}
				\dfrac{f(x+h)-f(x)}{h} & = &\displaystyle \dfrac{u(h)}{h} f'(x) + \sum_{m=1}^n\dfrac{v_m(h)}{h}w_m(x)
				\\
				f'(x) & = & \displaystyle\lim_{h\to 0}\left(\dfrac{u(h)}{h} f'(x) + \sum_{m=1}^n\dfrac{v_m(h)}{h}w_m(x)\right).
			\end{array}
		\end{equation*}
		Hence, if $f'\equiv0$, we need to have $f(x)=c$ for some constant $c$, and this case can be easily verified. Then, for $f'\not\equiv 0$, there exists $x\in\mathbb{C}$ such that $f'(x)\neq 0$, which implies $ \displaystyle\lim_{h\to 0}\left(\dfrac{u(h)}{h}\right) = 1$.
		
		Utilizing such limit, note that
		\begin{align*}%\label{Eq4.21}
			\displaystyle \lim_{h\to0} \, \left( \sum_{k=1}^{x}\hspace{-0.5cm}\rightarrow \frac{f(k+h) - f(k)}{h} \right) 
			& = 
			\lim_{h\to0} \, \left( \sum_{k=1}^{x}\hspace{-0.5cm}\rightarrow \frac{u(h) f'(k) + \sum_{m=1}^nv_m(h)w_m(k)}{h} \right)
			\\
			& =
			\lim_{h\to0} \, \left( \frac{u(h)}{h} \sum_{k=1}^{x}\hspace{-0.5cm}\rightarrow  f'(k)+\sum_{m=1}^n\dfrac{v_m(h)}{h}\sum_{k=1}^{x}\hspace{-0.5cm}\rightarrow w_m(k)\right) \\
			& =
			\lim_{h\to0} \, \left( \frac{u(h)}{h} \sum_{k=1}^{x}\hspace{-0.5cm}\rightarrow  f'(k) \right) + \lim_{h\to0} \, \left(  \sum_{m=1}^n\dfrac{v_m(h)}{h}\sum_{k=1}^{x}\hspace{-0.5cm}\rightarrow w_m(k)\right) \\
			& =
			\left(\sum_{k=1}^{x}\hspace{-0.5cm}\rightarrow  f'(k)\right) \cdot \lim_{h\to0} \, \left( \frac{u(h)}{h} \right) +\sum_{m=1}^n\left(\left( \sum_{k=1}^{x}\hspace{-0.5cm}\rightarrow w_m(k)\right) \cdot \lim_{h\to0} \,  \left( \frac{v_m(h)}{h} \right) \right)\\
			& =
			\sum_{k=1}^{x}\hspace{-0.5cm}\rightarrow  f'(k).
		\end{align*}
		Finally, by applying the derivative with respect to $x$ at $\displaystyle \sum^x\hspace{-0.45cm}\rightarrow f$, we obtain
		\begin{align*}%\label{Eq4.20}
			\displaystyle
			\frac{d}{dx} \Bigl( \sum^x\hspace{-0.45cm}\rightarrow f \Bigr)
			& =
			\displaystyle
			\lim_{h\to0} \, \frac{1}{h} \left( \sum_{k=1}^{x+h}\hspace{-0.5cm}\rightarrow f(k) - \sum_{k=1}^x\hspace{-0.5cm}\rightarrow f(k) \right) \\
			& \stackrel
			{A1}{=}
			\displaystyle
			\lim_{h\to0} \, \frac{1}{h} \left( \sum_{k=1}^{h}\hspace{-0.5cm}\rightarrow f(k)
			+ \sum_{k=1+h}^{x+h}\hspace{-0.65cm}\rightarrow f(k)
			- \sum_{k=1}^x\hspace{-0.5cm}\rightarrow f(k) \right)  \\
			& \stackrel{A2}{=}
			\displaystyle
			\lim_{h\to0} \, \frac{1}{h} \left( \sum_{k=1}^{h}\hspace{-0.5cm}\rightarrow f(k) + \sum_{k=1}^{x}\hspace{-0.5cm}\rightarrow f(k+h)
			- \sum_{k=1}^x\hspace{-0.5cm}\rightarrow f(k) \right) \\
			& =
			\displaystyle
			\lim_{h\to0} \, \left(\frac{1}{h} \sum_{k=1}^{h}\hspace{-0.5cm}\rightarrow f(k)+ \sum_{k=1}^{x}\hspace{-0.5cm}\rightarrow \frac{f(k+h) - f(k)}{h} \right) \\
			& \stackrel{A5^*}{=}
			\displaystyle
			\mathrm{ess}(f)
			+  \sum_{k=1}^{x}\hspace{-0.5cm}\rightarrow  f'(k) \, .
		\end{align*}
	\end{proof}
	
	\begin{corollary}\label{Corollary4.01} 
		Let $\mathcal{U}$be a summable set such that $\mathbb{C}[x]\subset \mathbf{S}(\mathcal{U})$. For $P\in\mathbb{C}[x]$, we have
		\begin{equation}\label{Eq4.19.1}
			\displaystyle
			\frac{d}{dx} \left( \sum^{x}\hspace{-0.45cm}\rightarrow P \right) = \mathrm{ess}(P) + \sum^{x}\hspace{-0.45cm}\rightarrow P' \, .
		\end{equation}
	\end{corollary}
	\begin{proof}
		For $P(x)=1$ the result is trivial. For $P(x)=x^a$ with $a\geq 1$, note that
		\begin{equation*}
			P(x+h)-P(x)=\underbrace{h}_{u(h)} P'(x)+\left(\sum_{m=1}^{a-2}\underbrace{h^{m+1}}_{v_m(h)}\underbrace{\binom{a}{m+1}x^{a-(m+1)}}_{w_m(x)}\right)+\underbrace{h^{a}}_{v_{a-1}(h)}\cdot \underbrace{1}_{w_{a-1}(x)}.
		\end{equation*}
		In case in which $P(x)$ is not a monomial, it suffices to apply linearity.
	\end{proof}
	
	In Theorem~\ref{tayloress}, we deal with derivative of fractional sums of holomorphic functions for which its successive derivatives satisfy the same hypothesis of Theorem~\ref{Theorem4.02} for a ``fixed'' summable set $\mathcal{U}$.
	
	\begin{theorem}\label{tayloress}
		Let $f$ be an holomorphic function. If for all $\lambda\geq 0$ there exist $2n+1$ functions $u_\lambda, v_{(1,\lambda)},\dots,v_{(n,\lambda)},w_{(1,\lambda)},\dots,w_{(n,\lambda)}:\mathbb{C}\to\mathbb{C}$ satisfying
		\begin{itemize}
			\item[i) ]$\forall\,m\in\{1,\dots,n\},\,\displaystyle\lim_{h\to0}\left(\dfrac{v_{(m,\lambda)}(h)}{h}\right)=0$,
			\item[ii) ] $\forall\,(x,h)\in\mathbb{C}\times\mathbb{C},\,\displaystyle f^{(\lambda)}(x+h)-f^{(\lambda)}(x)=u_\lambda(h)f^{(\lambda+1)}(x)+\sum_{m=1}^n v_{(m,\lambda)}(h)w_{(m,\lambda)}(x)$,
		\end{itemize}
		and if there is a summable set $\displaystyle\mathcal{U}$ such that $\displaystyle\bigcup_{\lambda\in\mathbb{N}}\left\{f^{(\lambda)},w_{(1,\lambda)},\dots,w_{(n,\lambda)}\right\}\subset \mathbf{S}(\mathcal{U})$,
		then, for all $y\in\mathbb{C}$, we have
		\begin{equation}\label{Eq4.28}
			\sum^{x}\hspace{-0.45cm}\rightarrow f=\sum_{k=1}^\infty{\left(\mathrm{ess}\left(f^{(k-1)}\right)+\displaystyle\sum^y\hspace{-0.45cm}\rightarrow f^{(k)}\right)}\dfrac{(x-y)^k}{k!}.
		\end{equation}
		In particular, for $y=0$, we obtain
		\begin{equation}\label{Eq4.29}
			\displaystyle 
			\sum^{x}\hspace{-0.45cm}\rightarrow f=\sum_{k=1}^\infty{\mathrm{ess}\left(f^{(k-1)}\right)}\dfrac{x^k}{k!}.
		\end{equation}
	\end{theorem}
	\begin{proof}
		It is sufficient to write the Taylor series of $\displaystyle\sum^{x}\hspace{-0.45cm}\rightarrow f$ by using Theorem~\ref{Theorem4.02}.
	\end{proof}
	
	\begin{corollary}\label{Corollary4.02} 
		Let $\mathcal{U}$ be a summable set such that $\mathbb{C}[x]\subset\mathbf{S}(\mathcal{U})$. For $P\in\mathbb{C}[x]$ of degree $n$, we have
		\begin{equation}\label{Eq4.19.2}
			\displaystyle
			\sum^{x}\hspace{-0.45cm}\rightarrow P = \sum_{k=1}^{n+1}\mathrm{ess}\left(P^{(k-1)}\right)\dfrac{x^k}{k!}.
		\end{equation}
	\end{corollary}
	\begin{proof}
		It suffices to adapt Corollary~\ref{Corollary4.01} to see that $(P,\mathcal{U})$ satisfies the hypotheses of Theorem~\ref{tayloress}.
	\end{proof}
	
	In Theorem~\ref{Theorem4.04}, we show that for $f(x) = \frac{1}{x^a}, a \in \mathbb{N}^*$, the derivative for it fractional sums can be evaluated, and the result agrees with Eq.~\eqref{Eq4.01}. Note that this is an example of function that is not covered by Theorems~\ref{Theorem4.02} and \ref{tayloress}.
	
	\begin{theorem}\label{Theorem4.04}
		Let $X=\{-1, -2, -3, \dots\}$. For any $a\in\mathbb{N}^*$, let $f_a(x)=\dfrac{1}{x^a}$. If there exists a summable set $\mathcal{U}$ such that $\{f_a,f_a'\}\subset\mathbf{S}(\mathcal{U})$, then $\mathrm{ess}(f_a)=a\zeta(a+1)$ and  
		\begin{equation}
			\dfrac{d}{dx}\sum^x\hspace{-0.45cm}\rightarrow f_a=\mathrm{ess}(f_a)+\sum^x\hspace{-0.45cm}\rightarrow f_a'
		\end{equation}
		for all $x\in\mathbb{C}-X$.
	\end{theorem}
	\begin{proof}
		By the fundamental formula for fractional summable functions, the fractional sums of $f_a$ and $f_a'$ are well defined (see Remark~\ref{somenotation}):
		\begin{equation*}
			\begin{array}{rcl}
				\displaystyle 
				\sum^x\hspace{-0.45cm}\rightarrow f_a 
				& = & 
				\begin{cases}
					\displaystyle\sum_{k=1}^\infty\left(\dfrac{1}{k^a}-\dfrac{1}{(k+x)^{a}}\right) & \text{ if } x\in\mathbb{C}-X,
					\\
					\displaystyle\dfrac{1}{(-x)^a}+\sum_{\substack{k=1 \\ k\neq -x}}^\infty\left(\dfrac{1}{k^a}-\dfrac{1}{(k+x)^{a}}\right) & \text{ if } x\in X,
				\end{cases}
				\\
				\text{and}&&
				\\
				\displaystyle \sum^x\hspace{-0.45cm}\rightarrow f_a' 
				& = & 
				\begin{cases}
					\displaystyle -a\sum_{k=1}^\infty\left(\dfrac{1}{k^{a+1}}-\dfrac{1}{(k+x)^{a+1}}\right) & \text{ if } x\in\mathbb{C}-X,
					\\
					\displaystyle \dfrac{-a}{(-x)^{a+1}}-a\sum_{\substack{k=1 \\ k\neq -x}}^\infty\left(\dfrac{1}{k^{a+1}}-\dfrac{1}{(k+x)^{a+1}}\right) & \text{ if } x\in X.
				\end{cases}
			\end{array}
		\end{equation*}
		Next, from \cite[p. 259-260]{Abramowitz1964}, %Eq.(6.3.16) and Eq. (6.4.10)
		it is clear that for all $x\in\mathbb{C}-X$, $\displaystyle \sum^x\hspace{-0.45cm}\rightarrow f_a=C_a+\dfrac{(-1)^{a+1}}{(a-1)!}\psi^{(a-1)}(x+1)$, where $\psi^{(0)}:=\psi$, $C_1=\gamma$, and $C_a=\zeta(a)$ for $a>1$. Therefore, for all $x\in\mathbb{C}-X$, we obtain
		\begin{equation*}
			\begin{array}{rcl}
				\displaystyle \dfrac{d}{dx}\sum^x\hspace{-0.45cm}\rightarrow f_a & = & \dfrac{d}{dx}\left(C_a+\dfrac{(-1)^{a+1}}{(a-1)!}\psi^{(a-1)}(x+1)\right)=\dfrac{(-1)^{a+1}}{(a-1)!}\psi^{(a)}(x+1)
				\\
				& = & \displaystyle\dfrac{(-1)^{a+1}}{(a-1)!}\left(\dfrac{a!}{(-1)^a}\left(-C_{a+1}+\sum^x\hspace{-0.45cm}\rightarrow f_{a+1}\right)\right)=-a\left(-\zeta(a+1)+ \sum^x\hspace{-0.45cm}\rightarrow \dfrac{f_a'}{-a}\right)
				\\
				& = & \displaystyle a\zeta(a+1)+\sum^x\hspace{-0.45cm}\rightarrow f_a'=\mathrm{ess}(f_a)+\sum^x\hspace{-0.45cm}\rightarrow f_a'.
			\end{array}
		\end{equation*}
	\end{proof}
	
	\begin{remark}
		To evaluate the essence of $f_a$ it is sufficient to see that  $f_a\in\mathbf{S}(\mathcal{U})$, since from \cite[p. 260]{Abramowitz1964}, %Eq.~(6.4.9)
		we have $\psi^{a}(0+1)=(-1)^{a+1}[a!\zeta(a+1)]$. This statement justifies writing $\mathrm{ess}(f_a)$ without specifying the summable set in question.
	\end{remark}
	
	\begin{theorem}\label{Theorem5}
		Let $f:\mathbb{C}\to\mathbb{C}$ and consider a summable set $\mathcal{U}$ such that $\{f,f'\}\subset\mathbf{S}(\mathcal{U})$. If $f$ satistisfies Eq.~\eqref{Eq4.19}, then for all $y\in\mathbb{C}$, the function $f_y\in\mathbf{T}(\mathcal{U})\subset\mathbf{S}(\mathcal{U})$, given by
		$f_y(x):=f(x+y)$, admits essence in $\mathcal{U}$, given by $\displaystyle\mathrm{ess}(f_y)=\mathrm{ess}(f)+\sum^y\hspace{-0.45cm}\rightarrow f'$. Moreover, $\displaystyle\dfrac{d}{dy}\left.\mathrm{ess}(f_y)\right|_{y=0}=\mathrm{ess}(f')$.
	\end{theorem}
	\begin{proof}
		\begin{equation*}
			\begin{array}{rcl}
				\mathrm{ess}(f_y) & = & \displaystyle\dfrac{d}{dx}\left.\left(\sum^x\hspace{-0.45cm}\rightarrow f_y\right)\right|_{x=0}
				\\
				& = & \displaystyle\dfrac{d}{dx}\left.\left(\sum_{k=1+y}^{x+y}\hspace{-0.65cm}\rightarrow f(k)\right)\right|_{x=0}
				\\
				& = & \displaystyle\dfrac{d}{dx}\left.\left(\sum^{x+y}\hspace{-0.5cm}\rightarrow f-\sum^{y}\hspace{-0.45cm}\rightarrow f\right)\right|_{x=0}
				\\
				& = & \displaystyle\left.\left(\mathrm{ess}(f)+\sum^{x+y}\hspace{-0.5cm}\rightarrow f'\right)\right|_{x=0}
				\\
				& = & \displaystyle\mathrm{ess}(f)+\sum^y\hspace{-0.45cm}\rightarrow f'.
			\end{array}
		\end{equation*}
		For the second part, just derive at $y=0$ in both sides and apply Proposition~\ref{Proposition4.01}.
	\end{proof}
	
	%
	% ----- Section 5  -----
	%
	\section{Examples and Applications}\label{Section5}
	
	In this section, we apply the results developed in Section~\ref{Section4} to give some examples and applications. Specifically, we: (i) provide one method to evaluating essences for functions $f$ for which the sum $\displaystyle \sum^{x} \hspace{-0.45cm}\rightarrow f$ is not known a priori; (ii) expand the well-known Euler-Maclaurin summation formula for a real/complex summation boundary limit; and (iii) propose a new method to evaluate divergent series.
	
	\subsection{Evaluating essences of functions from fractional sums theory}
	
	In what follows, we present a method for evaluating the essence of a function for which Proposition~\ref{Proposition4.02} cannot be applied directly. For a given function $\displaystyle F(x) = \sum^{x}\hspace{-0.45cm}\rightarrow f(k)$, in order to evaluate its essence we assume that $\mathrm{ess}(f) = u$, and proceed similarly as in Examples~\ref{Example5.05} and \ref{Example4.03}.
	
	\begin{example}\label{Example5.05}
		If $z\in\mathbb{C}$ is such that $e^z\neq 1$, then the function $f(x)= e^{z x}$ satisfies the hypotheses of Theorem~\ref{tayloress}.
		In fact, for $f^{(\lambda)}$, it suffices to take the functions $u_\lambda, v_{(1,\lambda)}, w_{(1,\lambda)}$, where $u_\lambda(h) = \dfrac{e^{z h}-1}{z}$ and $v_{(1,\lambda)}\equiv w_{(1,\lambda)}\equiv 0$. Then it follows that
		\begin{equation*}\label{Eq4.24}
			\displaystyle f^{(\lambda)}(x+h)-f^{(\lambda)}(x)= u_\lambda(h)f^{(\lambda+1)}(x)+\sum_{m=1}^1v_{(m,\lambda)}(h)w_{(m,\lambda)}(x).
		\end{equation*}
		Now, let $\mathcal{U}$ be a summable set such that $f\in\mathbf{S}(\mathcal{U})$, since $\displaystyle\bigcup_{\lambda\in\mathbb{N}}\left\{f^{(\lambda)},w_{(1,\lambda)}\right\}\subset\mathbf{S}(\{f\})\subset\mathbf{S}(\mathcal{U})$. Because $\mathrm{ess}(f^{(k)}) = z^{k}\mathrm{ess}(f)=z^ku$, from Theorem~\ref{tayloress} it follows that
		\begin{equation*}
			\begin{array}{rcl}
				\displaystyle \sum_{k=1}^{x}\hspace{-0.5cm}\rightarrow e^{z k}  & = &\displaystyle \sum_{k=1}^\infty\Bigl(z^{k-1}u\Bigr)\frac{x^k}{k!}
				\\
				& = & \dfrac{u}{z}\displaystyle \sum_{k=1}^\infty\frac{(z x)^k}{k!}
				\\
				& = & u\left(\dfrac{e^{z x}-1}{z}\right),
			\end{array}
		\end{equation*}
		and considering $x=1$, we get $\displaystyle e^{z} = u\Bigl(\frac{e^{z}-1}{z}\Bigr)$. Then
		$\mathrm{ess}(f)=\dfrac{z e^{z}}{e^{z} - 1}$ and we have 
		\begin{equation*}  
			\displaystyle
			\sum_{k=1}^{x}\hspace{-0.5cm}\rightarrow e^{z k}=\dfrac{e^{z}(e^{z x}-1)}{e^{z}-1}.
		\end{equation*}
		In particular, $\mathrm{ess}(x\mapsto e^{\pi i x})=\dfrac{\pi}{2}i$.  %This method will be used later to calculate the fractional sums of the $x^\alpha\sin(\beta x)$ and $x^\alpha\cos(\beta x)$. %Before do this, first we go to a generalization of the Theorem \ref{Theorem4.01}
	\end{example}
	
	\begin{remark}
		Note that functions of the type $x\mapsto u^x$, with $u\in\mathbb{C}-\{0,1\}$, are covered by Example~\ref{Example5.05}, since we only required that $f\in\mathbf{S}(\mathcal{U})$, and since $u^x:=e^{\ln(u)x}$.
	\end{remark}
	
	\begin{example}\label{Example4.03}
		Let $z\in\mathbb{C}$ such that $e^z\neq 1$; then $f(x) = e^{z x}x$ satisfies the hypotheses of Theorem~\ref{tayloress} for any summable set $\mathcal{U}$ such that $f\in\mathbf{S}(\mathcal{U})$. Indeed, note that if we consider the functions $T\in\mathbf{T}(\mathcal{U})\subset\mathbf{S}(\mathcal{U})$ defined by $T(x):=f(x+1)$, and $L\in\mathbf{L}(\mathcal{U})\subset\mathbf{S}(\mathcal{U})$ defined by $L(x):=e^{z}\cdot f(x)$, then $g(x):=e^{zx}$ can be written as a linear combination of $T$ and $L$: $g(x)=e^{-z}\cdot (T(x) - L(x))$. Hence, $g\in\mathbf{S}(\mathcal{U})$.
		
		Now, for $f^{(0)}$, just take $u_{0}(h)=\dfrac{e^{zh}-1}{z}$, $v_{(1,0)}(h)=he^{zh}-\dfrac{e^{zh}-1}{z}$, and $w_{(1,0)}(x)=e^{zx}$. Then it follows that $f$ satisfies the hypothesis (i) of Theorem~\ref{Theorem4.02}:
		\begin{equation}
			f^{(0)}(x+h)-f^{(0)}(x)=u_0(h)f^{(1)}(x)+\sum_{m=1}^1v_{(m,0)}(h)w_{(m,0)}(x).
		\end{equation}
		
		Note that, for $\lambda>0$, $f^{(\lambda)}$ can be expressed as a linear combination of $f$ and $e^{zx}$: $f^{(\lambda)}(x)=\lambda z^{\lambda - 1}e^{zx}+z^\lambda f(x)$, and both functions satisfy the hypothesis (i) of Theorem~\ref{Theorem4.02} when $\mathcal{U}$ is considered to be the summable set in this case. Therefore, $f$ satisfies Eq.~\eqref{Eq4.29}. Next, if we let $F(x)= \displaystyle \sum^{x}\hspace{-0.45cm}\rightarrow f$, $h(x) = e^{z x}$ and $u = \mathrm{ess}(f)$, then one has
		\begin{align*}%\label{Eq4.30}
			F(x) & = \displaystyle\sum_{k=1}^\infty \mathrm{ess}\left(f^{(k - 1)}\right)\dfrac{x^k}{k!}
			\\
			& = \displaystyle\sum_{k=1}^\infty \mathrm{ess}\left(z^{k - 1}\cdot f + (k - 1)z^{k-2}\cdot h\right)\dfrac{x^k}{k!}
			\\
			& = \displaystyle\sum_{k=1}^\infty \left(z^{k - 1}\cdot\mathrm{ess}(f)+(k - 1)z^{k - 2}\cdot\mathrm{ess}(h)\right)\dfrac{x^k}{k!}.
		\end{align*}
		Replacing $\mathrm{ess}(f)$ by $u$, we obtain
		\begin{align*}%\label{Eq4.31}
			F(x) & = \displaystyle\sum_{k=1}^\infty \left(z^{k - 1}u+(k - 1)z^{k - 2}\left(\dfrac{e^{z}}{e^{z} - 1}\lim_{h\to0}\left(\dfrac{e^{zh}-1}{h}\right)\right)\right)\dfrac{x^k}{k!}
			\\
			& = \displaystyle\dfrac{u}{z}\sum_{k=1}^{\infty}z^k\dfrac{x^k}{k!} + \dfrac{e^{z}}{e^{z} - 1}\sum_{k=1}^\infty \dfrac{(k - 1)z^{k - 1}x^k}{k!}
			\\
			& = \displaystyle\dfrac{u}{z}\left(-1 + \sum_{k=0}^{\infty}\dfrac{(z x)^k}{k!}\right) + \dfrac{e^{z}}{e^{z} - 1}\left(\sum_{k=1}^\infty \dfrac{z^{k - 1}x^k}{(k-1)!}-\sum_{k=1}^\infty \dfrac{z^{k - 1}x^k}{k!}\right)
			\\
			& = \displaystyle\dfrac{u}{z}\left(-1 + e^{z x}\right) + \dfrac{e^{z}}{e^{z} - 1}\left(x\sum_{k=0}^\infty \dfrac{(z x)^{k}}{k!}-\dfrac{1}{z}\left(-1+\sum_{k=0}^\infty \dfrac{(z x)^{k}}{k!}\right)\right)
			\\
			& = \displaystyle\dfrac{u}{z}\left(e^{z x} - 1\right) + \dfrac{e^{z}}{e^{z} - 1}\left(xe^{z x}-\dfrac{1}{z}\left(-1+e^{z x}\right)\right)
			\\
			& = \dfrac{u(e^{zx} - 1)}{z} + \dfrac{e^{z}(z x e^{z x}-e^{z x}+1)}{z (e^{z} - 1)}.
		\end{align*}
		Evaluating $F$ at $x=1$, we get
		\begin{equation*}%\label{Eq4.31}
			F(1) = \dfrac{u(e^{z} - 1)}{z} + \dfrac{e^{z}(z e^{z }-e^{z }+1)}{z (e^{z} - 1)};
		\end{equation*}
		and since $F(1) = e^{z}$, it follows that
		\begin{equation*}
			\begin{array}{rrcl}%\label{Eq4.32}
				& e^{z} & = & \dfrac{u(e^{z} - 1)}{z} + \dfrac{e^{z}(z e^{z }-e^{z }+1)}{z (e^{z} - 1)}
				\\
				\Rightarrow & \quad z(e^{z}-1)e^{z} -e^{z}(z e^{z}-e^{z}+1) &=& u(e^{z} - 1)^2
				\\
				\Rightarrow & e^{z}(z e^{z } - z - z e^{z}+e^{z}-1) & = & u(e^{z} - 1)^2
				\\
				\Rightarrow & u = \dfrac{e^{z}}{e^{z} - 1}\left(1-\dfrac{z}{e^{z} - 1}\right) .
			\end{array}
		\end{equation*}
		Therefore, $\mathrm{ess}(f)=\mathrm{ess}\left(x\mapsto e^{z x}x\right) = \dfrac{e^{z}}{e^{z} - 1}\left(1-\dfrac{z}{e^{z} - 1}\right)$. Moreover, since
		\begin{equation*}%\label{Eq4.33}
			\displaystyle
			F(x) = \dfrac{u(e^{zx} - 1)}{z} + \dfrac{e^{z}(z x e^{z x}-e^{z x}+1)}{z (e^{z} - 1)},
		\end{equation*}
		we conclude that
		\begin{equation*}%\label{Eq4.34}
			\begin{array}{rcl}
				\displaystyle
				\sum_{k=1}^x\hspace{-0.5cm}\rightarrow e^{z k}k &=& \left(\dfrac{e^{z}}{e^{z} - 1}\left(1-\dfrac{z}{e^{z} - 1}\right)\right)\dfrac{(e^{zx} - 1)}{z} + \dfrac{e^{z}(z x e^{z x}-e^{z x}+1)}{z (e^{z} - 1)}
				\\
				& = & \dfrac{e^{z}}{e^{z}-1}\left(xe^{z x}-\dfrac{e^{z x} - 1}{e^{z} - 1}\right).
			\end{array}
		\end{equation*}
		In particular, for $z = \pi i$, we have $\mathrm{ess}\left(x\mapsto e^{\pi i x}x\right)=\dfrac{1}{2}+\dfrac{\pi}{4}i$ and
		\begin{equation}%\label{Eq4.34}
			\displaystyle
			\sum_{k=1}^x\hspace{-0.5cm}\rightarrow e^{\pi i k}k = \frac{1}{2}e^{\pi i x}x + \frac{1}{4}e^{\pi i x} - \frac{1}{4}.
		\end{equation}
	\end{example}
	
	The results displayed in Propositions~\ref{Corollary5.01} and~\ref{Proposition5.01} can be obtained from Example~\ref{Example4.03}.
	
	\begin{proposition}\label{Corollary5.01}
		If $\mathcal{U}$ is a summable set such that $\{x\mapsto\cos(\pi x)x\}\subset \mathbf{S}(\mathcal{U})$, then it follows that
		\begin{equation*}
			\begin{array}{rcl}
				\displaystyle
				\sum_{k=1}^x\hspace{-0.5cm}\rightarrow \cos(\pi k)k &=& \displaystyle\frac{1}{2}\cos(\pi x)x + \frac{1}{4}\cos(\pi x) - \frac{1}{4} , \text{ and} 
				\\
				\displaystyle
				\sum_{k=1}^x\hspace{-0.5cm}\rightarrow \sin(\pi k)k& =& \displaystyle\frac{1}{2}\sin(\pi x)x + \frac{1}{4}\sin(\pi x) .
			\end{array}
		\end{equation*}
	\end{proposition}
	\begin{proof}
		Since $\sin(x)=\cos\left(x-\frac{\pi}{2}\right),\sin(-x)=-\sin(x)$ and $\cos(-x)=\cos(x)$, then, for all $z\in\mathbb{C}$ (in particular, for $z=\pi$), we have
		\begin{equation*}
			\begin{array}{rcl}
				e^{i(zx)}x&=&\cos(zx)x+i\underbrace{\cos\left(z\left(x-\dfrac{\pi}{2z}\right)\right)\left(x-\dfrac{\pi}{2z}\right)}_{\in\mathbf{T}(\mathbf{S}(\mathcal{U}))=\mathbf{S}(\mathcal{U})}+\dfrac{i\pi}{2z}\underbrace{\cos\left(z\left(x-\dfrac{\pi}{2z}\right)\right)}_{\in\mathbf{T}(\mathbf{S}(\mathcal{U}))=\mathbf{S}(\mathcal{U})} ,
				\\
				e^{i(-zx)}x&=&\cos(zx)x-i\cos\left(z\left(x-\dfrac{\pi}{2z}\right)\right)\left(x-\dfrac{\pi}{2z}\right)-\dfrac{i\pi}{2z}\cos\left(z\left(x-\dfrac{\pi}{2z}\right)\right).
			\end{array}
		\end{equation*}
		Therefore, $\{x\mapsto e^{izx}x,x\mapsto e^{-izx}x,x\mapsto e^{izx},x\mapsto e^{-izx}\}\subset\mathbf{S}(\mathcal{U})$. From the definitions $\sin(x)=\dfrac{e^{xi}-e^{-xi}}{2i},\cos(x)=\dfrac{e^{xi}+e^{-xi}}{2}$ and by linearity the result follows.
	\end{proof}
	
	\begin{proposition}\label{Proposition5.01}
		Let $x\in\mathbb{C}$ such that $e^\frac{2\pi i}{x}\neq 1$. If $\mathcal{U}$ is a summable set such that
		\begin{equation*}
			\left\{y\mapsto\cos\left(\dfrac{2\pi}{x}y\right)y\right\}\subset\mathbf{S}(\mathcal{U}),
		\end{equation*}
		then it follows that:
		\begin{equation}
			\displaystyle
			\sum_{k=1}^x\hspace{-0.5cm}\rightarrow \cos\left(\frac{2\pi k}{x}\right)k = \frac{x}{2} .
		\end{equation}
	\end{proposition}
	\begin{proof}
		Let $F(\beta)=\displaystyle
		\sum_{k=1}^\beta\hspace{-0.5cm}\rightarrow \cos\left(\frac{2\pi }{x}k\right)k$. Applying the definition of cosine, the linearity, the results from Example~\ref{Example4.03} and Corollary~\ref{Corollary4.02}, and setting $z=\frac{2\pi i}{x}$, we have
		\begin{equation*}
			\begin{array}{rcl}
				F(\beta)&=&\dfrac{1}{2}\left(\displaystyle
				\sum_{k=1}^\beta e^{zk}k+\displaystyle \sum_{k=1}^\beta e^{-zk}k\right)
				\\
				& = & \dfrac{1}{2}\left(\dfrac{e^{z}}{e^{z}-1}\left(\beta e^{z \beta}-\dfrac{e^{z \beta} - 1}{e^{z} - 1}\right) + \dfrac{e^{-z}}{e^{-z}-1}\left(\beta e^{-z \beta}-\dfrac{e^{-z \beta} - 1}{e^{-z} - 1}\right)\right)
				\\
				& = & \dfrac{1}{2}\left(\dfrac{e^{z}}{e^{z}-1}\left(\beta e^{z \beta}-\dfrac{e^{z \beta} - 1}{e^{z} - 1}\right) - \dfrac{1}{e^{z}-1}\left(\beta e^{-z \beta}-\dfrac{e^{-z \beta} - 1}{e^{-z} - 1}\right)\right).
			\end{array}
		\end{equation*}
		Putting $\beta = x$, we have $e^{\pm z\beta}=e^{\pm 2\pi i}=1$; hence
		\begin{equation*}
			\begin{array}{rcl}
				F(x)&=&\dfrac{1}{2}\left(\dfrac{e^{z}}{e^{z}-1}\left(x\cdot 1-\dfrac{1 - 1}{e^{z} - 1}\right) - \dfrac{1}{e^{z}-1}\left(x\cdot1-\dfrac{1 - 1}{e^{-z} - 1}\right)\right)
				\\
				& = & \dfrac{x}{2(e^z-1)}(e^z-1)
				\\
				& = & \dfrac{x}{2}.
			\end{array}
		\end{equation*}
	\end{proof}
	
	\begin{example}
		Let $\mathcal{U}$ be a summable set such that $\ln\in\mathbf{S}(\mathcal{U})$. Then $\mathrm{ess}(\ln)=-\gamma$. Moreover, if $\left\{\ln,x\mapsto\dfrac{1}{x},x\mapsto\dfrac{1}{x^2},x\mapsto\dfrac{1}{x^3},\dots\right\}\subset\mathbf{S}(\mathcal{U})$, then, for $|x|<1$,
		\begin{equation}\label{zetagamma}
			\ln\Gamma(x+1)=-\gamma x + \sum_{k=2}^\infty\dfrac{\zeta(k)}{k}(-x)^k.
		\end{equation}
	\end{example}
	\begin{proof}
		From \cite{Muller2011} it follows that for all $|x|<1$, $\displaystyle\ln\Gamma(x+1)=\sum^x\hspace{-0.45cm}\rightarrow \ln$. Therefore, by \cite[p. 258]{Abramowitz1964}, % Eq.~(6.3.2)
		$\mathrm{ess}(\ln)=\dfrac{d}{dx}\left.\ln\Gamma(x+1)\right|_{x=0}=\psi(0+1)=-\gamma$.
		
		To show \eqref{zetagamma}, note that from Theorem~\ref{Theorem4.04}, for $a\geq2$, we have:
		\begin{align*}
			\displaystyle
			\dfrac{d^a}{dx^a} \left[\ln\Gamma(x+1)\right]_{x=0} & = \dfrac{d^{a-1}}{dx^{a-1}} \left[\psi(x+1)\right]_{x=0}
			\\
			& = \dfrac{d^{a-1}}{dx^{a-1}}\left[-\gamma+\sum_{k=1}^x\hspace{-0.5cm}\rightarrow\dfrac{1}{k}\right]_{x=0}
			\\
			& =(-1)^{a}(a-2)!\cdot\mathrm{ess}\left(x\mapsto\dfrac{1}{x^{a-1}}\right)
			\\
			& = (-1)^{a}(a-1)!\zeta(a).
		\end{align*}
		Now, we write the Taylor series of $\ln\Gamma(x+1)$ at $0$.
	\end{proof}
	
	\begin{remark}
		Based on the previous example, it seems that the fractional summations have potential in discovering identities like \eqref{zetagamma}. % of this type.
	\end{remark}
	
	\begin{example}
		Let $f:\mathbb{C}\to\mathbb{C}$ given by $f(x)=\dfrac{1}{x(x+1)}$, if $x\notin \{0,-1\}$ and $f(x)=0$ otherwise. If $\mathcal{U}$ is a summable set such that $f\in\mathbf{S}(\mathcal{U})$, then
		$\displaystyle\mathrm{ess}(f)=1$ and for $x\neq -1,-2,-3,\dots$, \ $\displaystyle\sum^x\hspace{-0.45cm}\rightarrow f=\dfrac{x}{x+1}$.
	\end{example}
	\begin{proof}
		Note that $f$ is fractional summable:
		\begin{align*}
			\sum^x\hspace{-0.45cm}\rightarrow f& = \sum_{k=1}^\infty( f(k)-f(k+x) )
			\\
			&=\sum_{k=1}^\infty\left( \dfrac{1}{k(k+1)}-\dfrac{1}{(k+x)(k+x+1)}\right)
			\\
			&= \sum_{k=1}^\infty\left(\dfrac{1}{k}-\dfrac{1}{k+1}\right)-\sum_{k=1}^\infty\left(\dfrac{1}{k+x}-\dfrac{1}{k+x+1}\right)
			\\
			&=
			\lim_{n\to\infty}\left(\dfrac{1}{1}-\dfrac{1}{n+1}\right)-\lim_{n\to\infty}\left(\dfrac{1}{1+x}-\dfrac{1}{n+x+1}\right)
			\\
			& = 1-\dfrac{1}{1+x}
			\\
			&=\dfrac{x}{x+1}.
		\end{align*}
		Also, $\mathrm{ess}(f)=\dfrac{d}{dx}\left.\dfrac{x}{x+1}\right|_{x=0}=1$.
	\end{proof}
	
	\subsection{The Euler-Maclaurin summation formula for fractional sums}
	
	Here, we utilize the theory of fractional sums for obtaining generalizations of Euler-Maclaurin summation formula to a real or complex-value summation boundary. Theorem~\ref{Teo6} show below is a generalization of Euler-Maclaurin summation formula for polynomials.
	
	\begin{theorem}\label{Teo6}
		Let $\mathcal{U}$ be a summable set such that $\mathbb{C}[x]\subset\mathbf{S}(\mathcal{U})$ and $P\in\mathbb{C}[x]$. Then, for all $x\in\mathbb{C}$, we have
		\begin{equation}\label{Eulermac}
			\displaystyle
			\sum_{k=1}^{x}\hspace{-0.5cm}\rightarrow P(k) = \displaystyle\int_0^xP(t)\,dt +
			\frac{P(x)-P(0)}{2} + \sum_{m=1}^\infty\dfrac{B_{2m}}{(2m)!} \bigl(P^{(2m-1)}(x)-P^{(2m-1)}(0)\bigr).
		\end{equation}
	\end{theorem}
	\begin{proof}
		Suppose that $P$ has degree $n$ and let $Q\in \mathbb{C}[x]$ be a primitive of $P$, that is, $Q'=P$. %that is, $\displaystyle\int_0^xP(t)dt=Q(x)-Q(0)$.
		Therefore, from Corollary~\ref{Corollary4.02}:
		\begin{align}
			\displaystyle
			\sum^{x}\hspace{-0.45cm}\rightarrow P\hspace{0.05cm} & = \displaystyle \sum_{k=1}^{n+1} \mathrm{ess}\left(P^{(k-1)}\right)\dfrac{x^k}{k!} \notag
			\\
			& = \displaystyle \sum_{k=1}^{n+1} \mathrm{ess}\left(y\mapsto Q^{(k)}(y)\right)\dfrac{x^k}{k!} \notag
			\\
			%& =
			%\displaystyle\sum_{k=1}^{n+1}\mathrm{ess}\left(y\mapsto\sum_{m=0}^{\infty}\dfrac{Q^{(k+m)}(0)}{m!}y^m\right)\dfrac{x^k}{k!}
			%\\
			& = \displaystyle\sum_{k=1}^{n+1}\mathrm{ess}\left(y\mapsto\sum_{m=0}^{n}\dfrac{Q^{(k+m)}(0)}{m!}y^m\right)\dfrac{x^k}{k!} \label{before}
			\\
			& =
			\displaystyle\sum_{k=1}^{n+1}\left(\sum_{m=0}^{n}\dfrac{Q^{(k+m)}(0)}{m!}\mathrm{ess}\left(y\mapsto y^m\right)\right)\dfrac{x^k}{k!} \label{after}
			\\
			%& = \displaystyle\sum_{k=1}^{n+1} \sum_{m=0}^n\dfrac{B_m\cdot Q^{(k+m)}(0)x^k}{m!k!}
			%\\
			& = \displaystyle\sum_{m=0}^n\dfrac{B_m}{m!} \sum_{k=1}^{n+1}\dfrac{ Q^{(m+k)}(0)}{k!}x^k \notag
			\\
			%& = \displaystyle\sum_{m=0}^n\dfrac{B_m}{m!} \left(-Q^{(m)}(0)+\sum_{k=0}^{n+1}\dfrac{ Q^{(m+k)}(0)}{k!}x^k   \right)
			%\\
			& = \displaystyle\sum_{m=0}^n\dfrac{B_m}{m!} \left(Q^{(m)}(x)-Q^{(m)}(0)\right) \notag
			\\
			& = \displaystyle Q(x)-Q(0)+ \sum_{m=1}^\infty\dfrac{B_m}{m!}\bigl(P^{(m-1)}(x)-P^{(m-1)}(0)\bigr) \notag
			\\
			& = \displaystyle\int_0^xP(t)\,dt +
			\frac{P(x)-P(0)}{2} + \sum_{m=2}^\infty\dfrac{B_m}{m!} \bigl(P^{(m-1)}(x)-P^{(m-1)}(0)\bigr)\notag
			\\
			& =
			\displaystyle\int_0^xP(t)\,dt +
			\frac{P(x)-P(0)}{2} + \sum_{m=1}^\infty\dfrac{B_{2m}}{(2m)!} \bigl(P^{(2m-1)}(x)-P^{(2m-1)}(0)\bigr).\notag
		\end{align}
	\end{proof}
	
	Note that Eq.~\eqref{Eulermac} has exactly the same form already known to the Euler-Maclaurin summation formula for natural numbers $n$. This fact is interesting since we have now an interpretation for the left term.
	
	Eq.~\eqref{Eulermac} can be shown to be true for some functions, such as $\cos(zx)$ or $e^{zx}$, both with $e^z\neq 1$. However, for an arbitrary holomorphic function $P$, the reasoning from \eqref{before} to \eqref{after} becomes more complicated than that done with polynomials due to the fact that it is a series and not a finite sum.
	
	We illustrate this difficult here. We know that if $\displaystyle f_i(x)=1+\sum_{k=1}^i\dfrac{1}{k!}x^k$, then $f_i(x)$ converges pointwise to $e^x$ and $\displaystyle\sum_{k=1}^x\hspace{-0.5cm}\rightarrow f_i(k)$ also converges pointwise to $\displaystyle\sum_{k=1}^x\hspace{-0.5cm}\rightarrow e^k$. Therefore, one can conjecture that if a sequence of functions $f_i$ converges pointwise to a function $f$, then $\displaystyle\sum\hspace{-0.45cm}\rightarrow f_i$  also converges poitwise to $\displaystyle\sum\hspace{-0.45cm}\rightarrow f$. Unfortunately, this ``pointwise convergence conjecture'' is false (see open question (iv), Section~\ref{Section6}). For example, if $f(x)=0$ for all $x\in\mathbb{C}$ and for $i\in\mathbb{N}$ we define
	\begin{equation*}
		\begin{split}
			f_i:\mathbb{C}&\to\mathbb{C}
			\\
			x & \mapsto 
			\begin{cases}
				1, & \text{ if } x=i
				\\
				0, & \text{ elsewhere,}
			\end{cases}
		\end{split}
	\end{equation*}
	then, $f_i\to f$ pointwise, but for any $x\notin \mathbb{Z}$, we have
	\begin{equation*}
		\displaystyle\lim_{i\to\infty}\sum^x\hspace{-0.45cm}\rightarrow f_i=\lim_{i\to\infty}\sum_{k=1}^\infty f_i(k)-f_i(k+x)=\lim_{i\to\infty}1\neq 0 = \sum^x\hspace{-0.45cm}\rightarrow f.
	\end{equation*}
	
	\subsection{How to evaluate divergent series using fractional sums?}\label{DIVERGENT SERIES}
	
	An important issue related to divergent series is to be able to assign a unique value to it. The process is known as regularization of divergent series. There are several techniques to regularize a divergent series: the Cèsaro, Abel or Borel summation method, the smoothed sum method and the Ramanujan constant of a series, among others \cite{Chagas2021,Chagas2022}. In the following, we propose a way to regularize divergent series using the theory of fractional finite sums.
	
	Motivated by the following formal calculation for a fractional summable function $F$ asymptotically approximated by zero:
	\begin{align*}\label{Eq5.203}
		\mathrm{ess}(F)&=\lim_{h\to0}\dfrac{1}{h}\sum^h\hspace{-0.45cm}\rightarrow F
		\\
		&=\lim_{h\to0}\dfrac{1}{h}\sum_{k=1}^\infty\left(F(k)-F(k+h)\right)
		\\
		&=-\sum_{k=1}^\infty\left(\lim_{h\to0}\dfrac{F(k+h)-F(k)}{h}\right)
		\\
		&=-\sum_{k=1}^\infty F'(k),
	\end{align*}  
	that works for some functions, such as $F(x)=\dfrac{1}{x}$, we propose a method to regularize divergent series, as follows.
	Given a function $f$, if $\displaystyle \sum_{k=1}^{\infty} f(k)$ is divergent, taking a function $F:\mathbb{C}\to \mathbb{C}$ such that $F'=f$, we have
	\begin{equation}\label{Eq5.204}
		^\#\sum_{k=1}^\infty f(k):= -\mathrm{ess}(F).
	\end{equation}
	In other terms, we propose that
	\begin{equation}\label{Eq5.205}
		^\#\sum_{k=1}^\infty f(k) := - \Bigl( \frac{d}{dx} \sum^{x}\hspace{-0.45cm}\rightarrow F\,\Bigr) \Big|_{x=0}.
	\end{equation}
	
	Unfortunately, as well as occurs with other methods, the process of choice of the adequate function $F$ remains an open question. However, we found some examples of natural primitive functions $F$ for some functions $f$ for which the value assigned to the divergent series $\displaystyle \sum_{k=1}^{\infty} f(k)$ is known (see Examples~\ref{Example5.01}-\ref{Example5.04}).
	
	\begin{example}\label{Example5.01}
		If $f(x)=x$, we can take $F(x)=\dfrac{x^2}{2}$. Then,
		\begin{equation*}\label{Eq5.206}
			\mathrm{ess}(F) = \mathrm{ess}\left(\dfrac{x^2}{2}\right)=\dfrac{1}{2}\mathrm{ess}(x^2)=\dfrac{1}{2}\cdot\dfrac{1}{6}=\dfrac{1}{12},
		\end{equation*}
		that is,
		\begin{equation}\label{Eq5.207}
			^\#\sum_{k=1}^\infty k = -\dfrac{1}{12}.
		\end{equation}
	\end{example}
	
	\begin{example}\label{Example5.02}
		If $f(x)=\dfrac{1}{x}$, we can take $F=\ln$. Then,
		\begin{equation}\label{Eq5.209}
			\displaystyle
			\mathrm{ess}(F) = \mathrm{ess}(\ln)=-\gamma \quad \Rightarrow \quad ^\#\sum_{k=1}^{\infty} \frac{1}{k} = \gamma,
		\end{equation}
	\end{example}
	
	\begin{example}\label{Example5.03}
		If $f(x) = -e^{\pi i x}$, we can take $F(x)=-\dfrac{1}{\pi i} e^{\pi i x}$. Then, 
		\begin{equation*}\label{Eq5.211}
			\displaystyle
			\mathrm{ess}(F) = \mathrm{ess}\Bigl( - \frac{1}{\pi i}  e^{\pi i x} \Bigr) = -\frac{1}{\pi i} \cdot \mathrm{ess}((e^{\pi i})^x) = - \frac{1}{\pi i} \cdot \frac{e^{\pi i}\ln(e^{\pi i})}{e^{\pi i} - 1} = - \frac{1}{\pi i} \cdot \frac{\pi i}{2} = - \dfrac{1}{2}.
		\end{equation*}
		That is:
		\begin{equation}\label{Eq5.212}
			^\#\sum_{k=1}^\infty (-1)^{k+1} = \dfrac{1}{2},
		\end{equation}
	\end{example}
	
	\begin{example}\label{Example5.04}
		If $f(x)=-e^{\pi i x}x$, we can take $F(x) = -\dfrac{1}{\pi i}e^{\pi i x}x - \dfrac{1}{\pi^2}e^{\pi i x}$. Then, follows that
		\begin{align*}\label{Eq5.214}
			\displaystyle
			^\#\sum_{k=1}^\infty (-1)^{k+1}k
			= -\mathrm{ess} \left( -\dfrac{1}{\pi i}e^{\pi i x}x - \dfrac{1}{\pi^2}e^{\pi i x} \right)  = \dfrac{1}{\pi i}\mathrm{ess}(e^{\pi i x}x) + \dfrac{1}{\pi^2}\mathrm{ess}(e^{\pi i x}) = \dfrac{1}{4}.
		\end{align*}
	\end{example}
	
	\begin{remark}
		Note that if $F$ is a function such that $F'$ satisfies the Euler-Maclaurin summation formula \eqref{Eulermac} and $\displaystyle\sum_{k=1}^\infty F''(k)$ is divergent, then, formally, our regularization method gives us that
		\begin{equation*}
			\begin{array}{rcl}
				^\#\displaystyle\sum_{k=1}^\infty F''(k) & = & -\mathrm{ess}(F')
				\\
				& = & -\displaystyle\lim_{x\to0}\dfrac{1}{x}\left(F(x)-F(0)+\dfrac{F'(x)-F'(0)}{2}+\sum_{m=1}^\infty \dfrac{B_{2m}}{(2m)!}(F^{(2m)}(x)-F^{(2m)}(0))\right)
				\\
				& = & -\displaystyle F'(0) - \dfrac{F''(0)}{2} - \sum_{m=1}^\infty\dfrac{B_{2m}F^{(2m+1)}(0)}{(2m)!}
			\end{array}
		\end{equation*}
		In other words, such a method consists of simply ``erasing the terms that depend on infinity'':
		
		\begin{center}
			$\displaystyle\sum_{k=1}^\infty F''(k)$
			\\
			$\downarrow$
			\\
			$\displaystyle \cancel{F'(\infty)}-F'(0)+\dfrac{\cancel{F''(\infty)}-F''(0)}{2}+\sum_{m=1}^\infty \dfrac{B_{2m}}{(2m)!}\left(\cancel{F^{(2m+1)}(\infty)}-F^{(2m+1)}(0)\right)$
			\\
			$\downarrow$
			\\
			$ -\displaystyle F'(0) - \dfrac{F''(0)}{2} - \sum_{m=1}^\infty\dfrac{B_{2m}F^{(2m+1)}(0)}{(2m)!}$
			\\
			$\downarrow$
			\\
			$^\#\displaystyle\sum_{k=1}^\infty F''(k)$.
		\end{center}
	\end{remark}
	
	%
	% ----- Section 6 - Conclusions -----
	%
	\section{Final Remarks}\label{Section6}
	
	In this paper we have shown an open question proposed by Müller and Schleicher (see \cite{Muller2011}), concerning the existence of a constant obtained in derivative processes of a fractional sum. To do this task, we defined the concept of essence of a function. Moreover, we have proposed a method to compute the essence as well as fractional sums of functions by applying Taylor series. Additionally, with the theory of fractional summations, we can understand the limit boundaries of summation in the well-known Euler-Maclaurin summation formula that can be applied to real values. We have also proposed a new method for regularization of divergent series. To finalize the paper, we next present a list of open questions, in order to provide interesting questions to interested researchers.
	\begin{itemize}
		\item[(i) ]Is $\mathfrak{S}$ closed under finite union? That is, $X,Y\in\mathfrak{S}\Rightarrow X\cup Y\in\mathfrak{S}$?
		\item[(ii) ]Is $\mathfrak{S}$ closed under arbitrary union? That is, $S\subset\mathfrak{S}\Rightarrow \displaystyle\bigcup_{X\in S}X\in\mathfrak{S}$?
		\item[(iii) ]If the answer to (ii) is affirmative, then it is only necessary to deal with the largest summable set $\mathfrak{F}=\displaystyle\bigcup_{X\in \mathfrak{S}}X$, since the fractional sum relative to any other summable set $\mathcal{U}$ is a restriction of the ``largest fractional sum'' to the set $\mathbf{S}(\mathcal{U})\times\mathbb{C} \times\mathbb{C}$. Furthermore, which is that set?  
		\item[(iv) ] Is it possible to add hypotheses to the functions in the pointwise convergence conjecture in order to become true?
		\item[(v) ] Is it convenient to add new axioms to justify assertions such as: ``if $f_i\to f$ uniformly, then $\displaystyle\sum\hspace{-0.45cm}\rightarrow f_i\to \sum\hspace{-0.45cm}\rightarrow f$ uniformly/pointwise''?
	\end{itemize}
	
	%
	% ----- Section - Acknowledgements -----
	%
	\section*{Acknowledgements}
	
	This work was partially supported by the Brazilian agencies CNPq, CAPES and Funda\c{c}\~{a}o Arauc\'{a}ria $\#$PBA2022011000222.
	
	\section*{Conflict of interest}
	
	The authors declare that they have no conflict of interest.

\end{document}